\documentclass[sn-mathphys,Numbered]{sn-jnl}

\usepackage{graphicx}%
\usepackage{multirow}%
\usepackage{amsmath,amssymb,amsfonts}%
\usepackage{amsthm}%
\usepackage{mathrsfs}%
\usepackage[title]{appendix}%
\usepackage{xcolor}%
\usepackage{textcomp}%
\usepackage{xspace}
\usepackage{manyfoot}%
\usepackage{booktabs}%
\usepackage{algorithm}%
\usepackage{algorithmicx}%
\usepackage{algpseudocode}%
\usepackage{listings}%
\usepackage[capitalise]{cleveref}
\usepackage[disable]{todonotes}
\usepackage{siunitx}
\usepackage{a4wide}
\usepackage{subcaption}
\usepackage{comment}

\newcommand{\throatvar}{\Theta}
\newcommand{\calI}{\mathcal{I}}
\newcommand{\calN}{\mathcal{N}}
\newcommand{\calG}{\mathcal{G}}
\newcommand{\calE}{\mathcal{E}}

\newcommand{\dpthroat}{\Delta p_{c,ij}}
\newcommand{\dpthroatidx}[1]{\Delta p_{c,{#1}}}

\newcommand{\R}{\mathbb{R}}

\newcommand{\FIT}{FI-$\throatvar$\xspace}
\newcommand{\FIR}{FI-R\xspace}
\newcommand{\FIN}{FI-N\xspace}

\theoremstyle{plain}

\newtheorem{definition}{Definition}

\theoremstyle{remark}
\newtheorem{remark}{Remark}

\definecolor{darkgreen}{rgb}{0,0.5,0}

\DeclareSIUnit\bar{bar}
\begin{document}

\title[Article Title]{Improvement of fully-implicit two-phase pore-network models by employing generalized flux functions with additional throat variables.}

\author*[1]{\fnm{Martin} \sur{Schneider}}\email{martin.schneider@iws.uni-stuttgart.de}
\author[1]{\fnm{Hanchuan} \sur{Wu}}
\author[1]{\fnm{Maziar} \sur{Veyskarami}}
\author[2]{\fnm{Sorin} \sur{Pop}}
\author[1]{\fnm{Rainer} \sur{Helmig}}
\affil[1] {\orgdiv{Institute for Modelling Hydraulic and Environmental Systems}, \orgname{University of Stuttgart}, \orgaddress{\street{Pfaffenwaldring 61}, \city{Stuttgart}, \postcode{70569}, \country{Germany}}}
\affil[2] {\orgdiv{Faculty of Sciences}, \orgname{Hasselt University}, \orgaddress{\street{Agoralaan Gebouw D}, \city{Diepenbeek}, \postcode{3590}, \country{Belgium}}}

\abstract{In fully-implicit two-phase pore-network models, developing a well-converged scheme remains a major challenge, primarily due to the discontinuities in the phase conductivities. This paper addresses these numerical issues by proposing a generalized flux function that establishes a continuous flux expression for two-phase flows by introducing an additional throat variable $\throatvar$. Two approaches for expressing this additional throat variable are introduced: the first applies regularization strategies, while the second constructs an additional residual constraint equation. It is shown that this approach significantly improves accuracy and ensures the temporal convergence, as demonstrated through various numerical examples.}

\maketitle
\section{Introduction}
Pore-network models (PNMs) have always been considered as attractive tools to model and investigate transport in porous media, which account for pore-scale phenomena and provide a reasonable accuracy with low computational costs among approaches for pore-scale modeling \citep{blunt2001flow}. Unlike direct numerical simulation methods, which explicitly resolve pore-scale geometries, PNMs simplify void spaces into, for instance, interconnected pore bodies and pore throats \citep{raoof2013poreflow} or interconnected capillary tubes with varying radii \citep{chen2020fully, veyskarami2016modeling}, significantly reducing computational cost. PNMs can be broadly categorized as single-phase models and two/multi-phase models, which focus on the transport in single-phase and two/multi-phase flows respectively \citep{blunt2001flow}. Furthermore, two-phase PNMs can be classified into quasi-static and dynamic models \citep{celia1995recent}. Quasi-static PNMs describe the fluid distribution when the system has reached equilibrium, i.e., assume capillary pressure equilibrium and neglect viscous forces. Dynamic PNMs, however, take both capillary and viscous forces into account, allowing simulation of transient flow behavior, i.e., flow conditions between two equilibrium states \citep{koplik1985two}. Dynamic PNMs can be formulated using either single-pressure or two-pressure formulation \citep{thompson2002pore}. Another classification refers to the main time-integration schemes used in pore-network modeling: explicit, implicit or a combination of them, e.g., IMPES scheme \citep{joekar2010non}. 

Striving for a numerically stable PNM, with less/no limitation in time-step size, which makes the model especially suitable for coupled systems, leads to the development of fully-implicit PNMs \citep{Weishaupt2021, an2020transition}. In order to fulfill this desire, a numerical description of pore-scale phenomena, e.g., fluid phase displacement in the throat (``invasion'' or ``snap-off''), is crucial \citep{Wu2024}, as they could lead to sudden significant local changes, affecting the numerical behavior of the system. In fact, what determines the accuracy and efficiency of a PNM is how the pore-scale phenomena occurring at the throat are described. 

In the fully-implicit PNM proposed by \cite{weishaupt2020model,Weishaupt2021}, the decision of whether the wetting phase in the throat is displaced with the non-wetting phase or vice versa is made based on the solution of the last Newton iteration during the Newton loop of some time step. However, using the solution obtained from an iteration, which may not represent a fully converged solution, raises a concern about the reliability of such decision-making, affecting the accuracy of the model and even resulting in unphysical solutions at the end of the time step \citep{Wu2024}. In addition, that could lead to significant negative impacts on the convergence behavior, especially, when the estimations at the end of the iteration are still far away from a converged solution, as the decisions need to be changed repeatedly. Furthermore, since the change in the throat state could be reflected as a jump in the fluid phase fluxes from one iteration to the next one, it could make the convergence for the Newton solver very challenging and potentially leads to instability \citep{Wu2024}.

\cite{chen2020fully} introduced in their fully-implicit algorithm, which also used the solution estimated in each iteration for the decision-making, time-stepping strategies to ensure reaching a physical solution at the end of the time-step. Using this strategy, they estimated the proper time-step size for the next time step by using the fluxes from the previous time-step to keep the saturation in the physical range and allow only one throat to be invaded in each time-step. However, this predictor is only first-order accurate and may therefore fail to accurately predict the timing of invasions or snap-off events when using large time steps. 

In our previous paper, \cite{Wu2024}, we analyzed in detail the challenges and critical aspects related to fully-implicit PNMs and introduced regularization strategies for fluid flux estimation at the throat to improve accuracy and efficiency. There, we showed how determining the throat state based on the converged solution and not the previous iteration along with considering the previous throat state improves the accuracy in comparison to other methods.

In this study, we continue tackling the numerical challenges related to fully-implicit PNMs by proposing a generalized flux approximation strategy that allows for a continuous representation of flux at the throat in a two-phase pore-network model. Furthermore, we propose a novel algorithm that incorporates an additional variable, $\throatvar$, at the pore throats, which is simultaneously solved with the primary variables via a constructed residual, which greatly improves the accuracy and robustness of the fully-implicit PNM. Then, we analyze the temporal convergence of the proposed PNMs, discussing the related challenges and important aspects that need to be considered. 

This paper is structured as follows: in the next section, we comprehensively discuss the fundamentals of the fully-implicit PNMs, from network to mathematical model and pore-body and throat local relations. Then, we describe numerical aspects, e.g., time discretization, convergence behavior and flux approximation in \Cref{sec:numer}. We also present the new approach to estimate the fluxes using a new variable at the throat as well as flux regularization in the last part of this section. The results, analysis, and discussions are presented in \Cref{sec:results}. In the final section, we provide a summary and an outlook on the next possible steps in this field.

\section{Pore-network model}
\begin{figure}[h!]
    \centering
    \begin{subfigure}{0.45\textwidth}
        \centering
        \includegraphics[width=\textwidth]{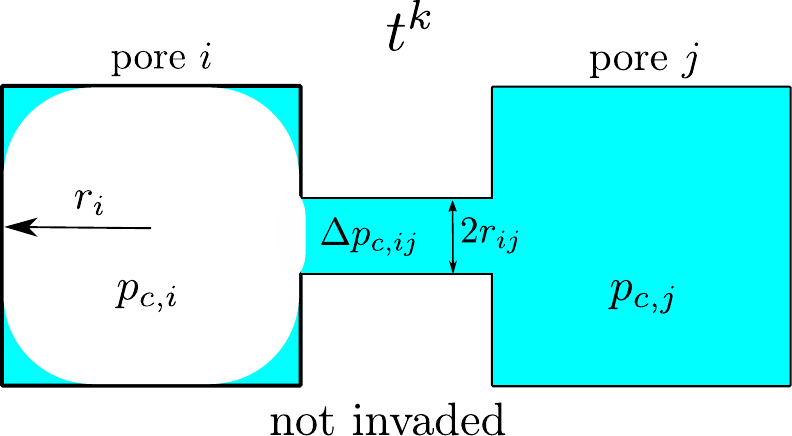}
        \caption{}
    \end{subfigure}
    \hfill
    \begin{subfigure}{0.45\textwidth}
        \centering
        \includegraphics[width=\textwidth]{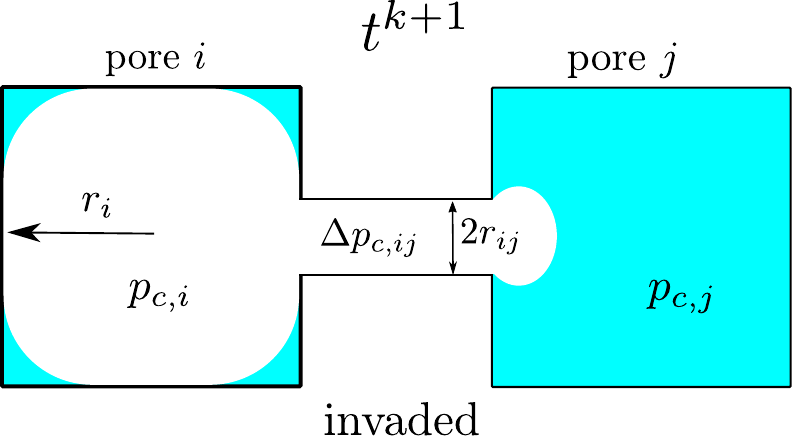}
        \caption{}
    \end{subfigure}
    \caption{\textbf{Schematic illustration of the pore-network model concepts.} The pore-body radius ($r_i$) and pore-throat radius ($r_{ij}$) are shown, assuming a cubic pore body and a circular throat cross-section and $\dpthroat = \max(p_{c,i},p_{c,j}) - p_{ce,ij}$, where $p_{ce,ij}$ is the entry capillary pressure: (a) when the pore throat is not invaded ($\dpthroat < 0$) at time $t^k$ and (b) when it is invaded ($\dpthroat \geq 0$) by the non-wetting phase at time $t^{k+1}$.}
    \label{fig:pore-throat}
\end{figure}

In this study, we focus on a pore network, that is a network of pore bodies, i.e., large void spaces, connected through pore throats. Such throats are defined as narrow void spaces, which are the primary source of resistance to flow.
\Cref{fig:pore-throat} illustrates pore bodies connected by a pore throat and the related parameters. In the following sections, we discuss different aspects of pore-network modeling in more detail. 

\subsection{Network} 
\label{subsec:network}
In the following, we introduce a mathematical definition of a network that consists of pore bodies connected by pore throats:
\begin{definition}[Network discretization]
\label{def:griddisc}
$\calG = (\calI, \calE)$ represents the network (graph), where
  \begin{enumerate}
  \item $\calI$ denotes the index set of pore bodies (vertices) with center $\mathbf{x}_i \in \R^3$ and volume $V_i$, $i \in \calI$.
  \item $\calE$ denotes the set of throats (edges) of the network with $e = \lbrace i,j \rbrace \in \calE$, $i,j \in \calI$. For simplicity, we often write $ij \in \calE$. 
  \item $\calN_i \subset \calI$, $i \in \calI$,  are the neighboring pore bodies connected to body $i$, i.e. $\calN_i:= \lbrace j \in \calI \; |  \; \exists e=\lbrace i,j \rbrace \in \calE \rbrace$.  
  \end{enumerate}
\end{definition}

In the following, quantities related to bodies are indicated by subscript $i\in \calI$, whereas those related to throats with subscript $ij \in \calE$.

\subsection{Mathematical model}
\label{subsec:model}
In this work, we focus on a dynamic fully-implicit PNM, describing the void space with pore bodies connected with pore throats, which uses a two-pressure formulation and a fully-implicit time-integration scheme. The pore bodies are treated as control volumes, at which the primary variables are located and the balance equations are solved. Such a network representation automatically yields a spatial discretization of the large pore spaces. 

The mass balance equation for each fluid phase $\alpha \in \lbrace w, n\rbrace$ within each pore body can be formulated as follows:

\begin{equation}
\label{eq:mass_balance}
    V_i \displaystyle\frac{ \mathrm{d} \left(  \rho_{\alpha, i} S_{\alpha,i}\right )}{\mathrm{d} t} + \sum_{j \in \calN_i}\left( \rho_\alpha Q_\alpha \right)_{ij} =V_i q_{\alpha,i}, \quad \forall i \in \calI, t \in (0, T],
\end{equation}
with phase saturation $S_\alpha$, density $\varrho_\alpha$, flux $Q_{\alpha}$, and source or sink term $q_\alpha$. Together with some appropriate initial conditions $S^0_\alpha, p^0_\alpha$ and the additional closure relations, \cref{eq:closures}, the system is closed.
\begin{equation}
    S_w + S_n = 1, \quad p_c = p_n-p_w,
    \label{eq:closures}
\end{equation}
where $p_c$ is the capillary pressure and $p_n$ and $p_w$ denote the non-wetting and wetting phase pressures, respectively.
As can be seen in \cref{eq:mass_balance}, the degrees of freedom, i.e., saturation and pressure, are assigned to the pore bodies. For the calculation of the volume flux $Q_{\alpha, ij}$ across the pore throat $ij \in \calE$ the following equation is used:
\begin{equation}
\label{eq:calculate_flux}
Q_{\alpha, ij} =  g_{\alpha, ij} \left( p_{\alpha, i}- p_{\alpha, j} \right), \; g_{\alpha, ij} = f \left( p_{c,i}, \ p_{c,j}\right),
\end{equation}
with throat conductivity $g_{\alpha, ij}$, given as a function of  $p_{c,i}$ and $p_{c,j}$ and depends on geometrical parameters, see \cref{subsec:conductivities}. Furthermore, the function $f$ is commutative, i.e. $f(x,y) = f(y,x)$ such that local flux conservation, $Q_{\alpha, ij} + Q_{\alpha, ji} = 0$, directly follows.
\Cref{{eq:mass_balance},{eq:closures},{eq:calculate_flux}} appear quite similar to those used in classical REV-scale models. However, despite their resemblance, the derivation of $p_c$ and $g_{\alpha, ij}$ is based on considerations of fluid interfaces, as explained in the following section.

\subsection{Pore-body and throat local relations}

\subsubsection{Capillary pressure--saturation relation}
\label{subsec:porerelations}
First, we need to mention, that in this paper, the displacement of the wetting phase by the non-wetting phase at the pore scale is called ``local drainage''. Similarly, the displacement of the non-wetting phase by the wetting phase at the pore scale is called ``local imbibition''.

The local capillary pressure of a pore body depends on the pore geometry and saturation. During the local drainage, when the non-wetting phase enters a pore body initially fully saturated with the wetting phase, the curvature of the phase interface is high, leading to a high capillary pressure. As wetting phase saturation decreases, the curvature decreases, causing a reduction in capillary pressure, which continues until the non-wetting phase contacts the pore-body walls and pushes the wetting phase into the corners. At this point, the capillary pressure starts to increase again, due to the further increase in the curvature. To simplify this non-monotonic $p_c-S_w$ relation, we use the approach introduced by \citep{joekar2010non}, where a monotonic relation is proposed for the local drainage process at the pore body $i$: 
\begin{equation}
\label{eq:pc-sw-vahid}
p_{c, i}(S_{w,i}) = \frac{2 \gamma}{r_{i} \left( 1 - \exp{\left( -6.83 S_{w,i} \right)} \right)}\, ,
\end{equation}
where $\gamma$ is the interfacial tension.

When the capillary pressure in the pore body reaches the entry capillary pressure of the connected pore throat, the non-wetting phase displaces the wetting phase in the pore throat, which we refer to as ``invasion''. Depending on the geometry of the pore throat, the wetting phase is fully displaced or only a small portion of it remains in the corners, leading to corner flow of the wetting phase. In other words, a pore throat is invaded, when it is fully filled with the non-wetting phase, e.g., in circular throats \citep{dias1986network, dias1986network2}, or both wetting and non-wetting phase exist in the angular throat \citep{al2005dynamic}, e.g., in square throats, with a higher conductivity of the non-wetting phase.

How the interface moves during the local imbibition in the pore body is described as the opposite of what occurs during the local drainage process. That means, \cref{eq:pc-sw-vahid} is sufficient to describe the $p_c-S_w$ relation of the pore body, when the wetting phase displaces the non-wetting phase. However, the local imbibition process in the pore throat is mainly influenced by two key mechanisms: cooperative pore filling and snap-off \citep{blunt2001flow}. Cooperative pore filling occurs when the pore-body-to-throat radius ratio is small, allowing the interface to bridge multiple throats in displacement of the non-wetting phase \citep{joekar2009simulating}. The ``snap-off'' occurs by an increase of wetting phase saturation resulting in a decrease in capillary pressure. The wetting phase (layers) on the wall and in the corners of the throats swell until reaching a state, when the wetting phase pinches the non-wetting phase off and fills the throat abruptly. The snap-off dominates the local imbibition process, when the pore-body-to-throat radius is large \citep{blunt2001flow}.

\subsubsection{Pore throat conductivities}
\label{subsec:conductivities}
One of the challenges, that two-phase pore-network models are facing with, is discontinuity of phase conductivities at the throats, when invasion or snap-off occurs. In these models, due to the rapid displacement processes in the throat, e.g., the Haines jump \citep{berg2013real, sun2019haines}, the interface movement through the throat, is not tracked.

\begin{figure}[h!]
	\centering
	\includegraphics[width=0.5\linewidth]{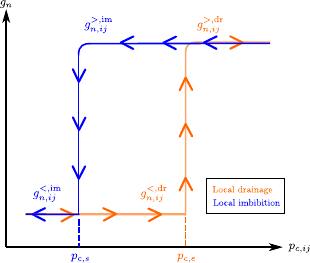}
	\caption{\textbf{Hysteresis in non-wetting phase conductivity.} The orange curve represents the non-wetting phase conductivity of pore throat during a throat local drainage process, where an invasion event occurs when $p_{c,ij} > p_{ce}$. In contrast, the blue curve depicts the conductivity during an imbibition process, where a snap-off event occurs when $p_{c,ij} < p_{c,s}$.}
	\label{fig:hysteresis}
\end{figure}

In local drainage, for instance, by overcoming the throat entry capillary pressure by pore-body capillary pressure due to increase in the non-wetting saturation, the non-wetting phase displaces the wetting phase in the throat, causing the non-wetting conductivity jumps from zero to a value near or equal to the maximum conductivity of the throat, depending on the throat cross-section. Simultaneously, the wetting phase conductivity drops from the maximum conductivity to near or equal to zero. Such an abrupt change of the meniscus due to the invasion is captured by the pore-network model with a discontinuous phase conductivity at the throat, see \cref{fig:hysteresis}.

In our model, the phase conductivity , $g_{\alpha,ij}$, of an invasion process is given by
\begin{equation}
\label{eq:conductivities}
    g_{\alpha,ij} = g_{\alpha,ij}^{>}H(\dpthroat)  + g_{\alpha,ij}^<(1-H(\dpthroat)), \quad \dpthroat = \max(p_{c,i},p_{c,j}) - p_{ce,ij}\, .
\end{equation}
The phase conductivities before and after the invasion are  $ g_{\alpha,ij}^<$, $g_{\alpha,ij}^>$, respectively. 
The Heaviside function,  $H(\dpthroat) := \mathbf{1}_{\dpthroat \geq 0}$, is equal to one after the invasion, where the function $\dpthroat$ (e.g. $\dpthroat = \max(p_{c,i},p_{c,j}) - p_{ce,ij}$) determines if the invasion occurs, dependent on pore-local capillary pressures. From the above definition, it becomes clear that $g_{\alpha,ij}$ has discontinuities at each time point $t_*$ where $\dpthroat$ changes its sign. 
Although \cref{eq:conductivities} is written for the invasion process, where we compare the throat capillary pressure to the entry pressure, a similar relation can be formulated for the snap-off process by comparing it with the snap-off capillary pressure. This could be indicated with an additional superscript, i.e., $\dpthroat^{\text{invasion}}$, $ \dpthroat^{\text{snap-off}}$. 

However, in this manuscript, for the sake of clarity and readability, we will only focus on the invasion process, proposing that the snap-off could be treated with a similar approach. In this work, we employ circular throat cross-sections, for which snap-off does not occur and we have $g_{n,ij}^< = 0$ for the non-wetting phase. Nevertheless, we will keep its contribution in the formulation, such that its generalization to other cross-sectional shapes remains straightforward. This is further discussed in \Cref{sec:GenFluxFunction}. For more details, we also refer to our previous work, \cite{Wu2024}.

\section{Numerical approximation}
\label{sec:numer}
Within this section, we discuss the temporal discretization of the ODE system in \cref{eq:mass_balance}.
Looking at the previous studies reveals the need for a comprehensive discussion on how to numerically handle the discontinuities in phase conductivities for a fully-implicit time discretization scheme. 
First, the implicit Euler discretization scheme and related issues are discussed. Then, we introduce a new concept of a generalized flux function to overcome these issues.

\subsection{Time discretization}
\label{subsec:time_disc}
The time interval $(0, T]$ is discretized into $N$ sub-intervals $(t^k,\,t^{k+1})$, where $k=0,\cdots,N-1$, with step sizes $\Delta t^{k+1} = t^{k+1} - t^k$. Furthermore, we define 
\begin{equation}
\label{eq:timepoints}
    \mathcal{T}_N := \lbrace t^k,  k=0,\cdots,N \rbrace
\end{equation}
as the set of time points.  
Integration of \cref{eq:mass_balance} over 
$(t^k,\, t^{k+1})$ and division by $\Delta t^{k+1} $ yield for each body $i \in \calI$:
\begin{equation}
\label{eq:timediscint}
    V_i \displaystyle\frac{ \left(  \rho_{\alpha, i} S_{\alpha,i}\right )^{k+1} - \left(  \rho_{\alpha, i} S_{\alpha,i}\right )^{k}}{ \Delta t^{k+1}} + \sum_{j \in \calN_i} \frac{1}{\Delta t^{k+1}} \int^{t^{k+1}}_{t^k} \left( \rho_\alpha Q_\alpha \right)_{ij} \, \mathrm{d}t = \frac{1}{\Delta t^{k+1}} \int^{t^{k+1}}_{t^k} V_i q_{\alpha,i} \, \mathrm{d}t.
\end{equation}
In the following, we denote with $(\bar\cdot)^{k+1}$ the time-averaging operator over the interval $(t^k,t^{k+1})$. 
Since $q_{\alpha,i}$ is a given function, the right-hand side can easily be calculated and is in the following denoted as $V_i \bar{q}^{k+1}_{\alpha,i}$.

For pore-network models, explicit or implicit Euler schemes have been extensively used to discretize the time-averaged flux integral. Due to the severe time-step size restrictions by using explicit schemes, the implicit or semi-implicit schemes are usually preferred, especially for diffusion-dominated (low capillary numbers) systems. As mentioned in \citep{gjennestad2018stable}, an efficiency increase of three orders of magnitude has been observed when using a semi-implicit Aker-type \citep{aker1998} pore-network model for cases with low capillary numbers compared to explicit schemes.

In this work, we focus on fully-implicit schemes, being unconditionally stable regarding the time-step size, while making the treatment of the flux discontinuities discussed in the previous section more challenging, as discussed in the following.

Applying the implicit Euler method on \cref{eq:timediscint}, results in:
\begin{equation}
\label{eq:implicitEuler}
    V_i \displaystyle\frac{ \left(  \rho_{\alpha, i} S_{\alpha,i}\right )^{k+1} - \left(  \rho_{\alpha, i} S_{\alpha,i}\right )^{k}}{ \Delta t^{k+1}} + \sum_{j \in \calN_i}\left( \rho_\alpha Q_\alpha \right)^{k+1}_{ij} = V_i \bar{q}^{k+1}_{\alpha,i}, \quad Q_{\alpha,ij}^{k+1} = g^{k+1}_{\alpha, ij} \left( p^{k+1}_{\alpha, i}- p^{k+1}_{\alpha, j} \right).
\end{equation}
As \cref{eq:implicitEuler} shows, the flux and, accordingly, $g^{k+1}_{\alpha, ij}$  implicitly depend on the solution of time level $t^{k+1}$. 
Using \cref{eq:conductivities} in \cref{eq:implicitEuler} gives
\begin{equation}
    Q_{\alpha,ij}^{k+1} = \left[g^{k+1,>}_{\alpha,ij}H(\dpthroat(t^{k+1}))  + g^{k+1, <}_{\alpha,ij}(1-H(\dpthroat(t^{k+1})))\right]\left( p^{k+1}_{\alpha, i}- p^{k+1}_{\alpha, j} \right),
    \label{eq:eulerfluxcasesHeaviside}
\end{equation}
or equivalently 
\begin{equation}
Q_{\alpha,ij}^{k+1} =
\label{eq:eulerfluxcases}
\begin{cases}
g^{k+1, <}_{\alpha,ij} \left( p^{k+1}_{\alpha, i}- p^{k+1}_{\alpha, j} \right), \quad &  \dpthroat^{k+1} <0, \\
g^{k+1,>}_{\alpha,ij} \left( p^{k+1}_{\alpha, i}- p^{k+1}_{\alpha, j} \right), \quad &  \dpthroat^{k+1} \geq 0.
\end{cases}
\end{equation}

The fluxes described by \cref{eq:eulerfluxcases} are discontinuous at each time point where invasion or snap-off occurs, i.e. the sign of $\dpthroat$ changes. This obviously leads to severe issues for derivative-based nonlinear solvers dealing with such systems. A simple strategy, often used in literature, to overcome these issues is to replace $g^{k+1}_{\alpha, ij}$ by $g^{k}_{\alpha, ij}$, i.e., the conductivity of the previous time step. With this, the fluxes no longer implicitly depend on conductivities such that discontinuities are prevented. However, especially for large or moderate time-step sizes, this might result in poor flux approximations and thus in inaccurate results. Additionally, when using this strategy, preventing unphysical saturations (negative or above one) cannot be guaranteed, which especially becomes important for physically complex flow processes (non-isothermal, compositional) where phase changes may occur. 

In the algorithms proposed by \cite{Weishaupt2021}, \cite{ an2020transition}, and \cite{chen2020fully}, the sign change of $\dpthroat$ within the Newton loop is allowed. However, in order to make Newton's method converge, only a single sign change is allowed such that the decision of being invaded or not is fixed and can't be corrected in later Newton iterations. This may lead to significant errors, as $\dpthroat$ might be completely overestimated within the first Newton iterations (unconverged results). Using this strategy in various test cases, we have observed large errors with unphysical behavior. Improvements can be achieved by including time-stepping strategies and repetition of time steps, or by applying large damping factors to achieve a monotonic increase of $\dpthroat$.

\begin{figure}[h!]
	\centering
	\includegraphics[width=0.7\linewidth]{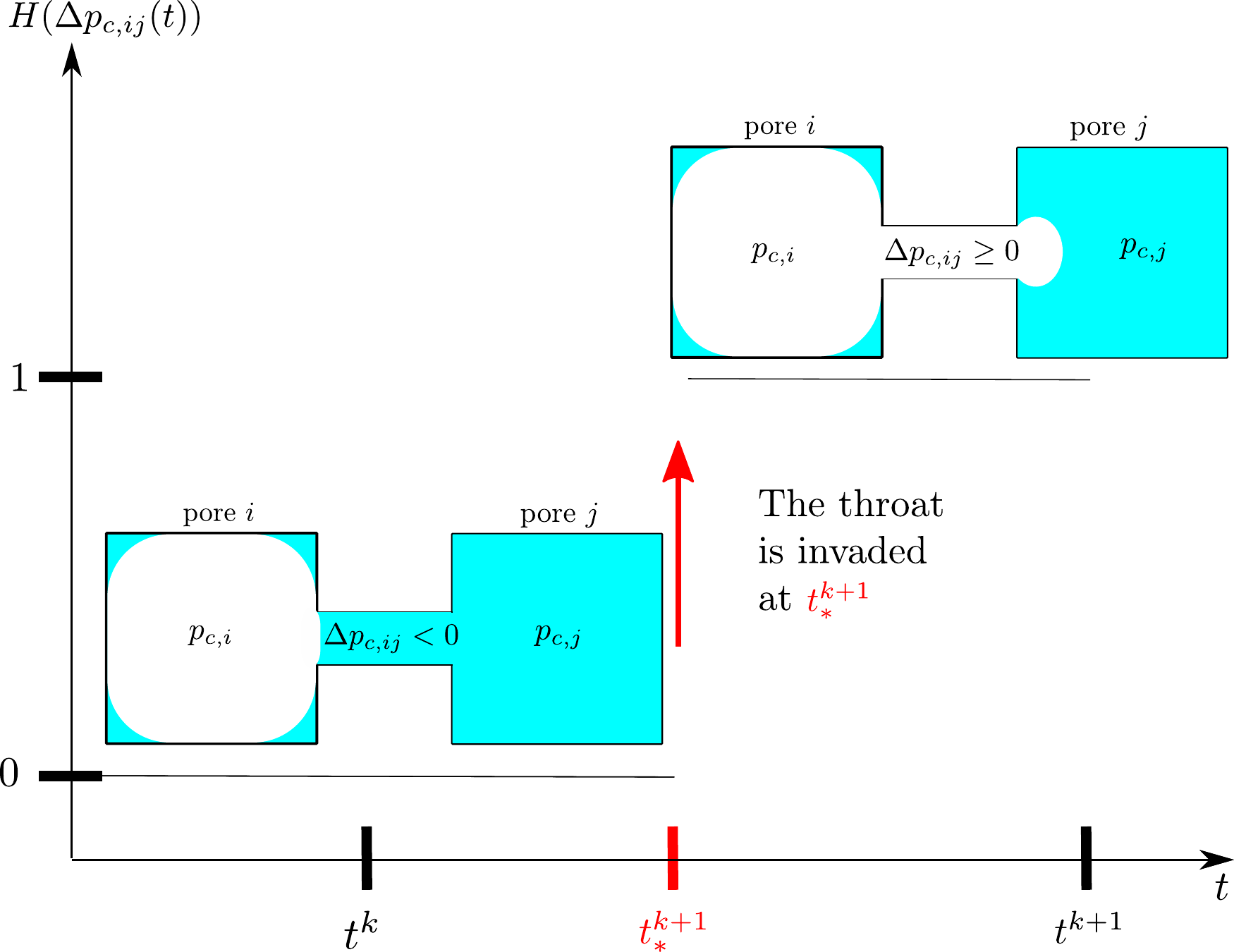}
	\caption{Schematic illustration of the throat invasion during the time interval $[t^k,t^{k+1}]$. The throat gets invaded at time $t^{k+1}_*$, at which $\dpthroat$ changes sign such that the Heaviside function transitions from zero (before invasion) to 1 (after invasion).}
	\label{fig:HvsDt}
\end{figure}

Another issue is related to the over- and underestimation of the flux during a time step, when using the flux approximation described in \cref{eq:eulerfluxcasesHeaviside}. Let us assume, for example, that an invasion event occurs within the time interval $(t^k,t^{k+1})$ in some throat $ij$, at time $t^{k+1}_*$ such that $t^k < t^{k+1}_* < t^{k+1}$ and stays invaded for the whole time step, which means $H(\dpthroat(t^k)) = 0$ but $H(\dpthroat(t^{k+1})) = 1$, equivalent to the second case in \cref{eq:eulerfluxcases}, see \cref{fig:HvsDt}. However, if $g^{>}_{\alpha,ij} >> g^{<}_{\alpha,ij}$, then the implicit Euler flux can lead to an overestimation of the time-averaged flux such that for the discrete solution we don't converge to $\dpthroat^{k+1} \geq 0$. However, since the Heaviside function has only two possible states, this would indicate that the decision of invasion, i.e., $H(\dpthroat(t^{k+1})) = 1$,  was wrong. On the other hand, assuming  $\dpthroat(t^{k+1})  < 0$ results in a numerical solution, where $\dpthroat^{k+1} \geq 0$, which again indicates that we should have evaluated the Heaviside function differently. This issue stems from the approximation of the exact time-averaged flux. Meaning that if we only allow two possible states, “invaded'' or “not invaded'', then there exists no discrete solution, when applying the implicit Euler flux \cref{eq:eulerfluxcasesHeaviside}. This leads to a non-converging behavior and emphasizes the necessity of modifications to allow for a more accurate flux approximation without the need to ``freeze'' the decision within the time step. In what follows, we show an example to elaborate on such a behavior.

\paragraph{Illustrative example}
As an example, we consider the following scenario: a local drainage process of two pore bodies connected by a throat, including two incompressible fluids as a part of a bigger network (see \cref{fig:HvsDt}).  The throat cross-section is assumed to be circular such that the non-wetting flux in the throat is blocked until $\dpthroat \geq 0 $ is reached, i.e. $g^{<}_{\mathrm{n},ij} = 0$. The non-wetting saturation increases in one of the pore bodies. Let us consider the two possible states described by \cref{eq:eulerfluxcases}. Assuming that no invasion happens, i.e. taking $g^{k+1, <}_{\alpha,ij}$ as conductivity (first state), might lead to a solution with $\dpthroat^{k+1} > 0 $, meaning that the invasion happened and we need to switch to the second case to calculate the fluxes. Then, by switching to the second case, we might end up with a solution with $\dpthroat^{k+1} < 0 $, which indicates that we need to switch back to the first case. Such a vicious circle is caused by an over- or underestimation of the exact flux and by allowing only two states during the time step. 

For this example,  $g^{<}_{\mathrm{n},ij} = 0$ and the time point at which $\dpthroat(t^{k+1}_*) = 0 $ is denoted by $t^{k+1}_*$, i.e., we assume that an invasion happens at throat $ij$ at time $t^{k+1}_*$ such that $t^k < t^{k+1}_* < t^{k+1}$ and that it stays invaded for the remaining of the time step, i.e. $\dpthroat(t) \geq 0$ for all $t^{k+1}_* \leq t \leq t^{k+1}$. Accordingly, the exact time-averaged non-wetting flux is given as:
\begin{equation}
\begin{aligned}
     \bar{Q}_{\mathrm{n},ij}^{k+1} &= \frac{1}{\Delta t^{k+1}} \int^{t^{k+1}}_{t^k}  Q_{\mathrm{n},ij}\, \mathrm{d}t = \frac{1}{\Delta t^{k+1}} \underbrace{\int^{t^{k+1}_*}_{t^{k}} Q_{\mathrm{n},ij} \, \mathrm{d}t}_{=0}  +  \frac{1}{\Delta t^{k+1}} \int^{t^{k+1}}_{t^{k+1}_*} Q_{\mathrm{n},ij} \, \mathrm{d}t \\
     &= \frac{1}{\Delta t^{k+1}} \int^{t^{k+1}}_{t^{k+1}_*} Q_{\mathrm{n},ij} \, \mathrm{d}t.
\end{aligned}
\end{equation}
Since $Q_{\mathrm{n},ij}(t)$ is a continuous function after the invasion event, i.e., for all $t>t^{k+1}_*$, it is observed that $\bar{Q}_{\mathrm{n},ij}^{k+1} \rightarrow 0$ if $t^{k+1}_* \rightarrow t^{k+1}$, meaning that the average non-wetting flux is very small if the invasion happens close to the end of the time step. \Cref{eq:eulerfluxcases} is not able to describe such a behavior. 

After the invasion event happened, the throat is ``open'' for the non-wetting phase such that ${Q}_{\mathrm{n},ij}(t) \approx {Q}_{\mathrm{n},ij}(t^{k+1})$ yields a good approximation for all $t \in (t^{k+1}_*, t^{k+1}]$. Therefore, the average flux for the non-wetting phase can be estimated by:
\begin{equation}
    \bar{Q}_{\mathrm{n},ij}^{k+1}\approx  \frac{t^{k+1} - t^{k+1}_*}{\Delta t^{k+1}} {Q}^{k+1}_{\mathrm{n},ij} \, ,
\end{equation}
which yields a better approximation of the exact flux compared to the fully-implicit Euler flux approximation described by \cref{eq:eulerfluxcases}. Hence, we introduce a new variable $\throatvar$ as $\throatvar := \frac{t^{k+1} - t^{k+1}_*}{\Delta t^{k+1}}$, that can be interpreted as a flux scaling factor. $\throatvar$ tends to zero ($\throatvar \rightarrow 0$) if the invasion event occurs closely before the next time point $t^{k+1} $ and approaches one ($\throatvar \rightarrow 1$) if the invasion event happens shortly after $t^k$. 

This example motivates the idea of introducing a generalized flux approximation function, that depends on an additional variable, $\throatvar$, at the throats. Without this additional scaling factor, the exact time-averaged fluxes are strongly over- or underestimated for each time step, where invasion or snap-off occurs. 

\subsection{Generalized flux functions}
\label{sec:GenFluxFunction}
As highlighted in the previous section, we can interpret $\throatvar \in [0,1]$ as a scaling factor (indicator), which is zero, when $H(\dpthroat(t)) = 0$ and equal to one if $H(\dpthroat(t)) = 1$ within the whole time step. 

Proceeding similar as in the previous section and by using \cref{eq:eulerfluxcasesHeaviside}, allows to approximate the exact time-averaged flux as
\begin{equation}
\label{eq:implicitFlux}
\begin{aligned}
 \frac{1}{\Delta t^{k+1}} \int^{t^{k+1}}_{t^k} \left( \rho_\alpha Q_\alpha \right)_{ij} \, \mathrm{d}t  &\approx \frac{1}{\Delta t^{k+1}} \left[\int^{t^{k+1}}_{t^k} g_{\alpha, ij} \, \mathrm{d}t \right] \rho^{k+1}_{\alpha,ij}  \left( p^{k+1}_{\alpha, i}- p^{k+1}_{\alpha, j} \right), \\
 \frac{1}{\Delta t^{k+1}} \left[\int^{t^{k+1}}_{t^k} g_{\alpha, ij} \, \mathrm{d}t \right] & \approx g_{\alpha,ij}^{>,k+1} \bar{H}  + g_{\alpha,ij}^{<,k+1}(1- \bar{H}) \, ,
\end{aligned}
\end{equation}
where
\begin{equation}
\label{eq:AvgHeaviside}
\begin{aligned}
 \bar{H} &:= \frac{1}{\Delta t^{k+1}} \int^{t^{k+1}}_{t^k} H(\dpthroat) \, \mathrm{d}t \, .
\end{aligned}
\end{equation}
Since $H(\dpthroat) = 0$ before the invasion,$t^k < t < t^{k+1}_*$, and $H(\dpthroat) = 1$ after the invasion, $t^{k+1}_* < t < t^{k+1}$, \cref{eq:AvgHeaviside} gives, in fact, $\throatvar^{k+1}_{ij}$ and we can write:
\begin{equation}
\label{eq:thetainterpretation}
    \throatvar^{k+1}_{ij} = \frac{1}{\Delta t^{k+1}} \int^{t^{k+1}}_{t^k} H(\dpthroat) \, \mathrm{d}{t}.
\end{equation}
Therefore, the generalized implicit  flux function can be written in terms of $\throatvar_{ij}$ as
\begin{equation}
    \label{eq:generalizedFlux}
    Q^{k+1}_{\alpha,ij}(p_{\alpha,i},p_{\alpha,j}, \throatvar_{ij}) :=  \left[g_{\alpha,ij}^{>,k+1}\throatvar_{ij}^{k+1} + g_{\alpha,ij}^{<,k+1}(1-\throatvar_{ij}^{k+1}) \right] \left( p^{k+1}_{\alpha, i}- p^{k+1}_{\alpha, j} \right).
\end{equation}
By considering $\throatvar_{ij}$ as an additional variable, the flux approximation \cref{eq:generalizedFlux} is continuous with respect to the extended variable space, in contrast to the implicit Euler flux given by \cref{eq:eulerfluxcases}.

Using this generalized flux approximation, \cref{eq:generalizedFlux}, in \cref{eq:timediscint} results in the following novel fully-implicit discrete formulation:
\begin{equation}
\label{eq:timedisc}
    V_i \displaystyle\frac{ \left(  \rho_{\alpha, i} S_{\alpha,i}\right )^{k+1} - \left(  \rho_{\alpha, i} S_{\alpha,i}\right )^{k}}{ \Delta t^{k+1}} + \sum_{j \in \calN_i}  \rho^{k+1}_{\alpha,ij} Q^{k+1}_{\alpha,ij}(p_{\alpha,i},p_{\alpha,j}, \throatvar_{ij})  = V_i \bar{q}^{k+1}_{\alpha,i}.
\end{equation}

It is important to note that:
\begin{remark}
The closer $\throatvar$ to zero or one, the better the invasion/snap-off event is captured by the time-stepping strategy and, therefore, the smaller the flux approximation error. 
\end{remark}

\begin{remark}
    \label{rem:hysteresis}
    \Cref{eq:thetainterpretation,eq:generalizedFlux} are derived assuming weather only invasion or only snap-off occur within a sing time step. In case of occurring both invasion and snap-off within a time step, the time-step size needs to be refined, such that the assumption holds. 
\end{remark}

In the following, we will discuss how to express $\throatvar$ as an additional variable in \cref{eq:timedisc}. 
In \Cref{subsec:fluxreg}, we will consider the case, where $\throatvar$ is approximated by making use of \cref{eq:thetainterpretation}. Afterward, we will discuss the case where $\throatvar$ is introduced as an additional variable at the throat in \Cref{subsec:mdtheta}.

\subsubsection{Flux regularization}
\label{subsec:fluxreg}
One of the approaches to guarantee the smoothness of $\throatvar$ is to regularize the Heaviside function, which will be denoted by $H_\delta$ in the following. Replacing the Heaviside function with the regularized one in \cref{{eq:thetainterpretation}}, and considering $t$ as a continuous variable, instead of integration over a given time interval, the following relation can be derived:
\begin{equation}
\label{eq:thetaregcont}
    \throatvar_{ij}(t) = H_\delta(\dpthroat(t)).
\end{equation}

This gives a regularized flux approximation function, which is quite similar to the work we have recently presented in \cite{Wu2024}, although we employ regularized Heaviside functions here instead of applying monotone splines in the conductivity laws. 
Here, it is constructed such that $\throatvar_{ij}$ is only larger than zero if $\dpthroat \geq 0$. Thus, $H_\delta$ is derived by shifting a regularized Heaviside function. 

A simple expression for $H_\delta \in C^1(\R)$, is given by, cf. \cite{Akkerman2017}:
\begin{equation}
H_\delta(s) = 
\begin{cases}
0, \quad & s \leq 0, \\
\frac{1}{2}\left( 1 - \cos\left(\frac{\pi s}{\delta} \right) \right), \quad &  0 < s < \delta, \\
1, \quad  \phantom{-\delta < }\, &  \delta \leq s  \, .
\end{cases}
\end{equation}
When discretizing \cref{eq:thetaregcont}, the discrete values $\throatvar^{k+1}_{ij}$ can be interpreted as an approximation of the time-averaged integral in \cref{eq:thetainterpretation}. Here, we simply choose:
\begin{equation}
\label{eq:thetaregdisc}
    \throatvar^{k+1}_{ij}= H_\delta(\dpthroat^{k+1}).
\end{equation}

Doing so, \cref{eq:timedisc} gives the following regularized fully-implicit discrete formulation:
\begin{equation}
\label{eq:timedisc_reg}
    V_i \displaystyle\frac{ \left(  \rho_{\alpha, i} S_{\alpha,i}\right )^{k+1} - \left(  \rho_{\alpha, i} S_{\alpha,i}\right )^{k}}{ \Delta t^{k+1}} + \sum_{j \in \calN_i}  \rho^{k+1}_{\alpha,ij} Q^{k+1}_{\alpha,ij}(p_{\alpha,i},p_{\alpha,j}, H_\delta(\dpthroat))  = V_i \bar{q}^{k+1}_{\alpha,i}.
\end{equation}



\subsubsection{Additional throat variable}
\label{subsec:mdtheta}
In this section, we present a new approach, where $\throatvar$ is introduced as an additional variable for each throat $ij \in \calE$. Thus, we do not need to replace $\throatvar$ in \cref{eq:generalizedFlux}, but instead, construct an additional equation for each $\throatvar_{ij}$ to determine the throat invasion state. This additional equation is generally written, using a residual function, $R$, as follows:
\begin{equation}
    R_\throatvar(\dpthroat, \throatvar_{ij}) = 0.
\end{equation}

It is important to construct this residual function $R_\throatvar$ such that for $\dpthroat = 0$ it becomes zero for any $\throatvar_{ij}$, i.e.,  $R_\throatvar(0, \throatvar_{ij}) \equiv 0$. Thus, whenever an invasion event happens, the condition $\dpthroat = 0$ determines the value of $\throatvar_{ij}$ which must satisfy $0 < \throatvar_{ij}< 1$. 
Additionally, $R_\throatvar$ must be chosen such that $R_\throatvar(\lbrace \dpthroat < 0 \rbrace, 0) \equiv 0$ and  $R_\throatvar(\lbrace \dpthroat > 0\rbrace, 1) \equiv 0$ such that it captures the correct states of ``not invaded'' $\throatvar = 0$ and ``invaded'' $\throatvar = 1$. 


A simple residual function that captures these conditions is given by:
\begin{equation}
    R_\throatvar := (1-\throatvar_{ij})\max(0, \dpthroat) - \throatvar_{ij} \min(0, \dpthroat).
\end{equation}
Accordingly, the new discrete model is defined as:
\begin{subequations}
\label{eq:timedisc_md}
\begin{align}
    \label{eq:timedisc_mass_md}
    V_i \displaystyle\frac{ \left(  \rho_{\alpha, i} S_{\alpha,i}\right )^{k+1} - \left(  \rho_{\alpha, i} S_{\alpha,i}\right )^{k}}{ \Delta t^{k+1}} + \sum_{j \in \calN_i}  \rho^{k+1}_{\alpha,ij} Q^{k+1}_{\alpha,ij}(p_{\alpha,i},p_{\alpha,j}, \throatvar_{ij})  &= V_i \bar{q}^{k+1}_{\alpha,i}, \quad &&\forall \, i \in \calI \, ,\\
    (1-\throatvar^{k+1}_{ij})\max(0, \dpthroat^{k+1}) - \throatvar^{k+1}_{ij} \min(0, \dpthroat^{k+1}) &= 0, \quad &&\forall \, ij \in \calE.
    \label{eq:timedisc_throat_md}
    \end{align}
\end{subequations}
Here, the mass balance equation, \cref{eq:timedisc_mass_md}, is related to the pore bodies, whereas the additional equation, \cref{eq:timedisc_throat_md}, for $\throatvar$ is related to the throats. The system of equations presented by \cref{eq:timedisc_md} is then coupled. 

It is important to mention that \cref{eq:timedisc_throat_md} is only valid until the invasion occurs at the throat. After the invasion, \cref{eq:timedisc_throat_md} is replaced by another residual function that accounts for snap-off, which depends on $\dpthroat^\text{snap-off}$, the throat geometry and the snap-off phenomena. As mentioned before, in this work, we focus on the local drainage process and a network consisting of throats with circular cross-sections. Thus, the snap-off process is not discussed here and after the invasion the throat remains invaded, i.e. $\throatvar = 1$.

From \cref{eq:timedisc_throat_md} the following states can be derived:
\begin{subequations}
\begin{align}
    &\dpthroat^{k+1}  < 0  &\Leftrightarrow \; \; \; \;&\throatvar_{ij} = 0, &&\text{ i.e., no invasion within time step,} \\
    &\dpthroat^{k+1}   > 0 &\Leftrightarrow \; \; \; \; &\throatvar_{ij} = 1, &&\text{ i.e., invasion at } t^k, \label{eq:state_invadedbeginning}\\
    &\dpthroat^{k+1}   = 0  &\Leftrightarrow \; \; \; \; &\throatvar_{ij} \in (0,1), &&\text{ i.e., invasion at some } t^{k+1}_* \in (t^k,t^{k+1}). \label{eq:state_invaded}
\end{align}
\end{subequations}
\cref{eq:state_invadedbeginning} means that even though the throat is assumed to be invaded directly at the beginning of the time step, i.e. at time $t^k$, the capillary pressure in the throat obtained at the end of the time step is above the entry pressure and not equal to the entry capillary pressure of the throat. This behavior might be related, e.g., to the initial or boundary conditions but also by small pressure differences between neighboring pore bodies. This will be further discussed in \cref{sec:results}. 
\cref{eq:state_invaded} is the typical case for an invasion event, where at the time $t^k$ the throat is not yet invaded, but the invasion occurs at $t^{k+1}_* \in (t^k,t^{k+1})$. With this new approach, it is guaranteed that  $\dpthroat^{k+1}   = 0$ such that the throat entry capillary pressure is perfectly matched. The value of $\throatvar_{ij}$ provides an estimation of when the invasion happens, as discussed in \cref{subsec:time_disc}.

It should be pointed out that:
\begin{remark}
This new approach comes with the cost of having an additional variable at throats. Although, within each time step, most of the throats $\throatvar_{ij}$ yield trivial states $\lbrace 0,1\rbrace$ such that they could also be eliminated within each Newton iteration before solving the linearized system of equations, this is not done here and the additional variables are kept for the entire simulation. 
\end{remark}

\section{Numerical results}
\label{sec:results}
In this section, we focus on three test cases constructed for the verification of the presented pore-network models. Here, we specifically focus on a comprehensive analysis of the temporal convergence behavior. To the best of our knowledge, such an analysis has not been carried out before for the type of PNMs considered within this work. 
It is self-evident that a reliable numerical model must exhibit convergent behavior. However, there is lack of mathematical theory proving the temporal convergence and the related stability of pore-network models with discontinuous time events (invasion or snap-off).

To the best of our knowledge, the only work that discusses temporal convergence of PNMs is \cite{gjennestad2018stable}, which focuses on a different type of PNM called Aker-type model \citep{aker1998}. In their model, all volume is assigned to throats, and interfaces are explicitly tracked within them, assuming a continuous displacement process. This contrasts with the models presented in this study, where wetting phase displacement in the pore throats begins when the invasion criterion is fulfilled, resulting in a no-displacement state before and a full-displacement state after the invasion for a circular throat. As a result, their models do not exhibit the conductivity discontinuities considered in this study.
Additionally, \cite{gjennestad2018stable} examined temporal convergence only for a one-dimensional test case. In a second test case, which involves a single Haines jump within a 2D network, temporal convergence is not further discussed.

The main challenge to ensure convergence stems from the fact that the overall error is strongly influenced by the accurate prediction of invasion or snap-off events, which are modeled as discontinuities in conductivity in classical PNMs.
These events are characterized by rapid interface movement and have a non-local impact \citep{armstrong2013}. 
Thus, accurate prediction of these events is essential for achieving temporal convergence. Using explicit schemes with small time-step sizes, such issues can be circumvented. However, as mentioned before and also discussed by \cite{gjennestad2018stable}, explicit schemes are subject to strong time-step size restrictions for low capillary number cases, rendering them impractical for many applications.  

To enable a detailed and focused discussion of convergence behavior, the following test cases are kept simple. However, to the best of our knowledge, the second and third test cases exhibit more complexity than any other test cases considered in the literature regarding temporal convergence, since analysis is performed for networks in which multiple invasion events occur. The first test case is constructed such that the analytical solution is known. In the second test case, a heterogeneous bifurcating network is investigated. In the third test case, a fully connected, two-dimensional heterogeneous network is considered, with a specific focus on discussing the influence of event prediction accuracy on the overall $L^2$-error.

For all test cases, we consider cubic pore bodies and throats with circular cross-sections. By doing so, we focus on invasion events and eliminate the need to consider snap-off events. Other pore and throat shapes can also be directly incorporated into the presented models. However, as they do not offer additional insights for the analysis and discussion, they are not considered further here.

The error is analyzed using the following discrete $L^2$-error for variable $u = (S_{w,i})_{i\in \calI} $:
\begin{equation}
    E_u := \left( \frac{1}{|\calI|} \sum_{k} \Delta t^k \lVert u^k - u_\text{ref}(t^k) \rVert_2^2   \right)^{\frac{1}{2}}\, ,
\end{equation}
where $u_\text{ref}$ is the reference solution, either given by the analytical solution (if available) or by the solution calculated for small time-step sizes. In addition, for the reference solution, time steps are further refined just before and after invasion events to ensure that errors induced by these events are controlled. This also guarantees that the time point of invasion is matched within a given tolerance.
$u_\text{ref}(t^k)$ is then obtained by interpolation between reference time steps. 
$u^k = (u_i^k)_{i\in \calI}$ is a vector containing all discrete solution values at time $t^k$. The time-step size used on refinement level $i$ is indicated by a superscript if necessary.

For quantifying the error at time points of invasions, we define the set of time points and throats where invasions happen:
\begin{equation}
    \mathcal{T}^*_N := \lbrace \lbrace t^{k},ij \rbrace \mid \dpthroat^{k-1} < 0 \text{ and } \dpthroat^{k} \geq 0  \rbrace.
\end{equation}
With this, the error related to invasion events can be quantified as follows:
\begin{equation}
\label{eq:errorinvasions}
    E^*_{p_{ce}} := \left( \frac{1}{|\mathcal{T}^*_N|}  \sum_{\lbrace t^{k},ij \rbrace \in \mathcal{T}^*_N } |  \dpthroat^{k} |^2  \right)^{\frac{1}{2}}\, .
\end{equation}
Errors related to this are often referred to as \emph{prediction errors} in the following.

The following labels are used for the considered schemes:
\begin{itemize}
    \item \FIN: Scheme that switches the throat invasion state within Newton iterations \citep{weishaupt2020model}. 
    \item \FIR: Scheme that employs regularization with interval $\delta$, \cref{subsec:fluxreg},
    \item \FIT: Scheme with additional throat variable, \cref{subsec:mdtheta}.
\end{itemize}


\begin{remark}
\label{re:newtontimesteps}
It is important to note that the following analysis focuses on temporal convergence, i.e. decreasing $\overline{\Delta t}$. However, this does not imply that a constant $\Delta t$ is used throughout the simulations. Instead, a targeted time-step size is specified, and if convergence is not achieved within a given number of solver steps, Newton's method repeats the time steps with a refined time-step size. Therefore, smaller $\overline{\Delta t}$ do not necessarily indicate smaller local time steps, e.g. before or after invasion events. Thus, there is no direct link between $\overline{\Delta t}$ and prediction errors related to invasion events.
\end{remark}

\subsection{One-dimensional case}
\label{sec:onedimensional}
\todo[inline]{Maybe show same plot for setup as for the other cases.}
We consider a one-dimensional pore network consisting of 10 cubic pore bodies with the same inscribed radius of $200\ \si{\micro\metre}$. These pore bodies are connected by throats, each having an inscribed radius of $100\ \si{\micro\metre}$. Initially, the network is fully saturated with water. At the inlet of the network, a non-wetting phase fluid is injected at a constant rate of $5 \times 10^{-10} \ \si{\kilogram\per\second}$, and the drainage process takes $5 \times 10^2 \  \si{\second}$. At the outlet, $p_w = 10^{5}\ \si{\pascal}$ together with outflow conditions for the non-wetting phase are set.

\begin{figure}[h!]
    \centering
    \includegraphics[width=0.8\textwidth]{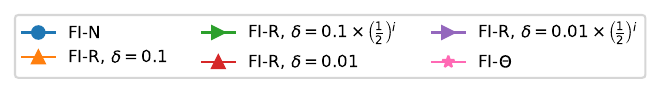}
    \vspace{0.1em}
    \begin{subfigure}{0.32\textwidth}
        \centering
        \includegraphics[width=\textwidth]{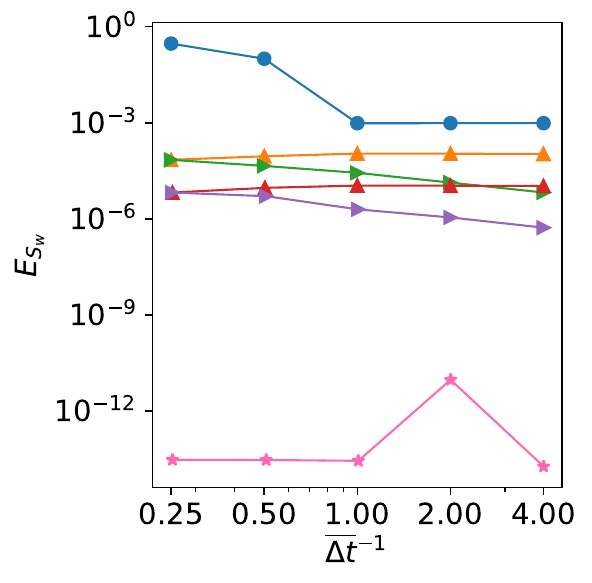}
    \end{subfigure}
    \begin{subfigure}{0.32\textwidth}
        \centering
        \includegraphics[width=\textwidth]{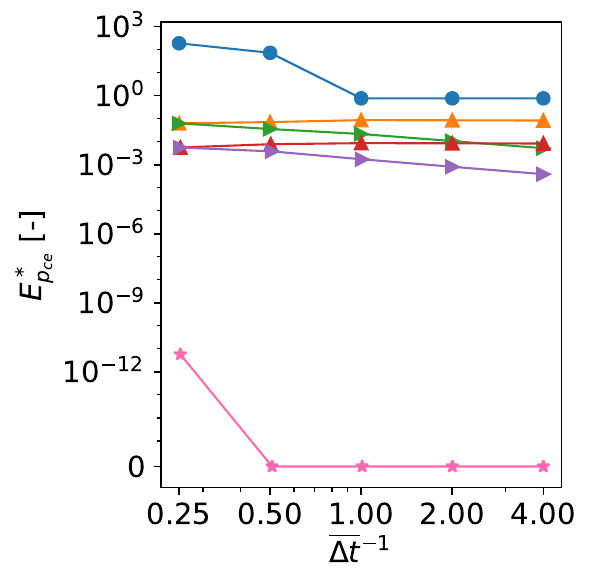}
    \end{subfigure}
    \hfill
    \begin{subfigure}{0.32\textwidth}
        \centering
        \includegraphics[width=\textwidth]{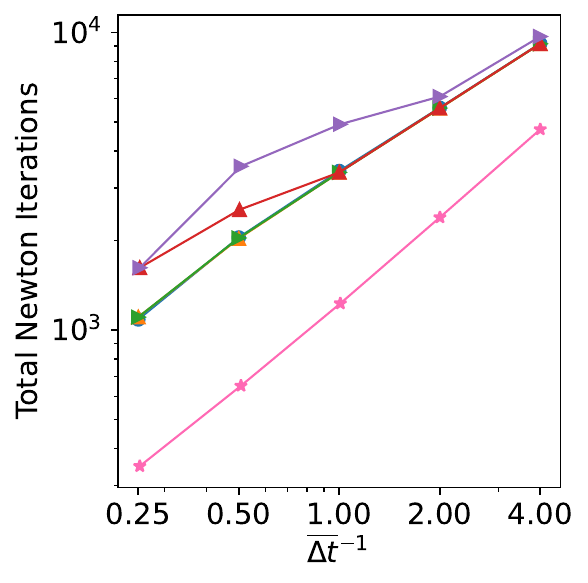}
    \end{subfigure}
    \caption{\textbf{One-dimensional test case.} $L^2$-errors (left), prediction errors of invasion events \cref{eq:errorinvasions} (middle), and total number of Newton iterations (right) plotted over the average time-step sizes. }
    \label{fig:1d_time_convergence}
\end{figure}

Firstly, we measure the accuracy of the simulation results with respect to the saturation and efficiency of each method throughout the simulation process. Accuracy is assessed by calculating the temporal $L^2$-error, between the numerical and analytical solutions for saturation, and the errors related to invasion events \cref{eq:errorinvasions}. Efficiency is evaluated by the total number of Newton iterations required for the whole simulation. The results are shown in \cref{fig:1d_time_convergence}. 

Looking at the results, we can see that the new \FIT scheme is capable of reproducing the exact solution (up to the given solver tolerance). This is expected when a constant influx is given since the \FIT scheme is locally mass conservative and it exactly captures the entry capillary pressure for invasion events \cref{eq:state_invaded}, as the middle plot of \cref{fig:1d_time_convergence} shows.
It can be seen in the left plot of \cref{fig:1d_time_convergence} that the \FIT scheme needs fewer Newton iterations compared to the other schemes.
In this test case, to investigate the convergence behavior of Newton's method, we set the maximum number of allowed Newton iterations to 500. This is to prevent the solver from repeating the time step with a reduced time-step size if the convergence requires a large number of Newton iterations, as it would normally be the case as detailed in \cref{re:newtontimesteps}. Therefore, \cref{fig:1d_time_convergence} clearly shows the better Newton convergence behavior for \FIT compared to the other schemes for this test case. 
\todo[inline]{Can we comment on that... does it also need less time steps? Then it is clear why it needs less Newton iterations. === I looked at  your data and it seems that it is not because of more time steps but because of better Newton convergence. The question would be if the worse convergence is caused by line search, which is not used by \FIT compared to the other schemes}

The errors of the \FIR schemes are approximately three orders of magnitude smaller than those of \FIN for small time-step sizes. For larger time steps, the errors of \FIN become significantly larger, whereas the errors remain stable for the \FIR schemes. 
Decreasing the regularization parameter for smaller time-step sizes in the \FIR schemes also leads to a decrease in $L^2$-errors. This is due to the correlation between $\delta$ and prediction accuracy, meaning that for a small $\delta$, it is expected that the deviation from entry capillary pressure $p_{ce}$ to be smaller. This relationship can be clearly seen when comparing the left plot with the middle plot of \cref{fig:1d_time_convergence}, where the reduction of prediction errors results in smaller $L^2$-errors.
Interestingly, it is also observed that the \FIR schemes do not exhibit a clear convergence behavior when $\delta$ is not decreased. This is because, for this simple one-dimensional test case, the overall error is primarily dominated by the accuracy in predicting the invasion events, i.e. $ |  \dpthroat^{k} |$ for all $\lbrace t^k,ij \rbrace \in \mathcal{T}^*_N$, see \cref{fig:1d_time_convergence} (middle). 
It can also be seen that the decreasing of $\overline{\Delta t}$ does not necessarily lead to smaller prediction errors. In general, this may be due to the Newton method locally refining time steps (see \cref{re:newtontimesteps}), though this is unlikely here since we allow a large number of iterations. Here, it is caused by invasion events occurring shortly before or directly after some time $t^k$, such that a simple refinement strategy does not effectively reduce prediction errors.

\subsection{Bifurcating network case}
In this section, we analyze two-dimensional bifurcating pore networks, which are designed as follows: all pores have the same size, while pore throats are assigned random radii, as shown in \cref{{fig:set_up_bifuracted}}. 
\begin{figure}[h!]
	\centering
	\includegraphics[width=0.8\linewidth]{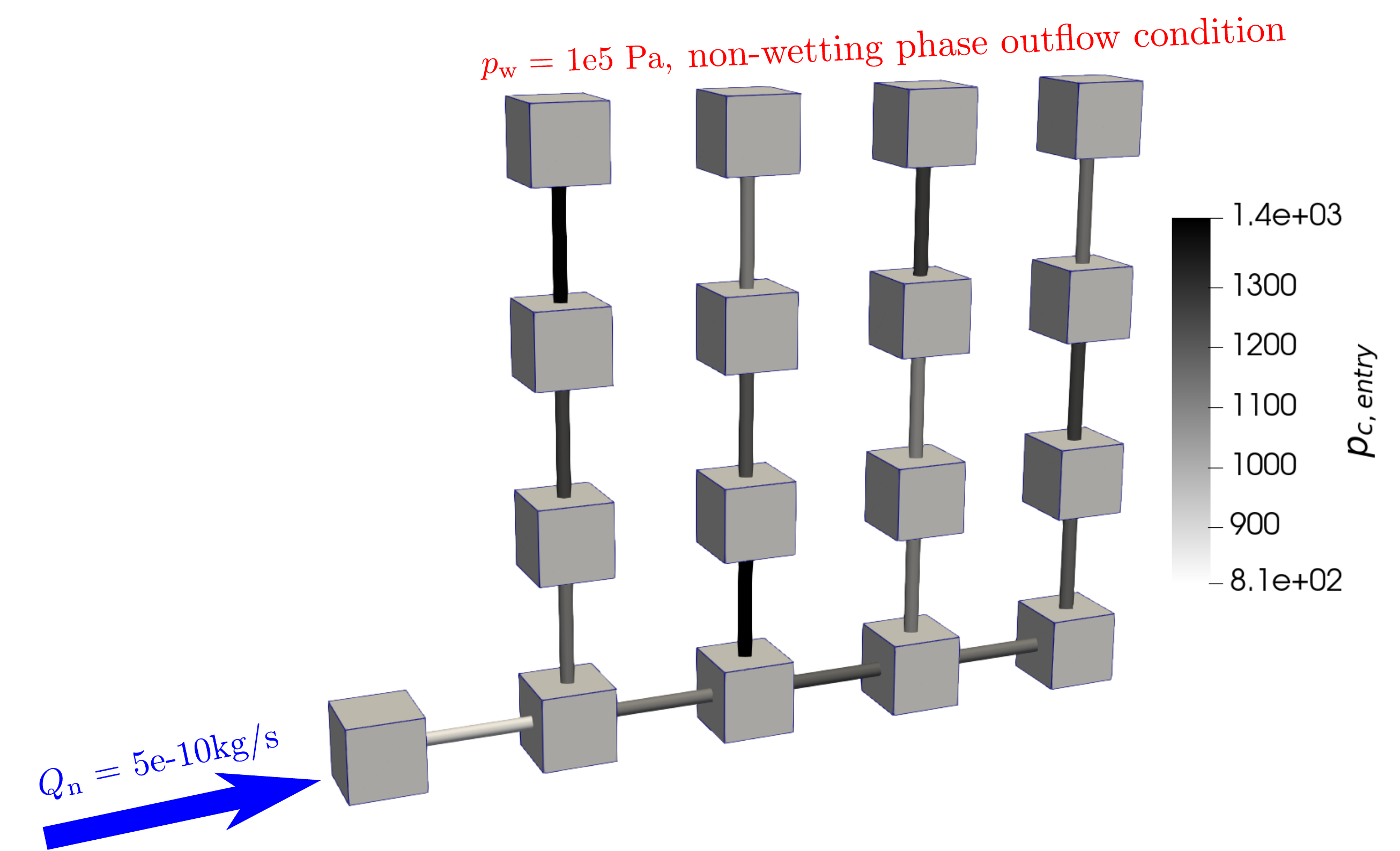}
	\caption{\textbf{Setup for the bifurcated case.} On the left bottom, non-wetting phase is injected with a fixed rate of $Q_n = $\SI{5e-10}{\kilogram\per\second}. At the four outlet pores on the top, a pressure of $p_\mathrm{w} = $ \SI{1e5}{\pascal} is set and the non-wetting fluid can flow out of the domain.}
	\label{fig:set_up_bifuracted}
\end{figure}
This design enables us to track the sequence of invaded throats. Because the next throat to be invaded is always the not invaded one with the largest radius, i.e., the lowest capillary entry pressure, that is connected to an already invaded pore body. Since the correct invasion sequence is known in advance, this case can be seen as a solid benchmark case. 

Initially, the entire network is saturated with water, and a non-wetting phase liquid is injected at the inlet at a fixed flow rate. We compare the accuracy of the different methods, focusing on whether they can correctly determine the invasion sequence and, if so, how accurate the numerically calculated capillary pressure matches the entry capillary pressures of the invaded throat in each invasion event.

\begin{figure}[h!]
    \centering
    \includegraphics[width=0.8\textwidth]{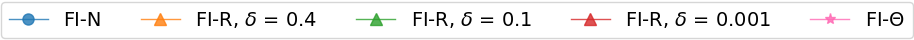}
    \begin{subfigure}{0.48\textwidth}
        \centering
        \includegraphics[width=\textwidth]{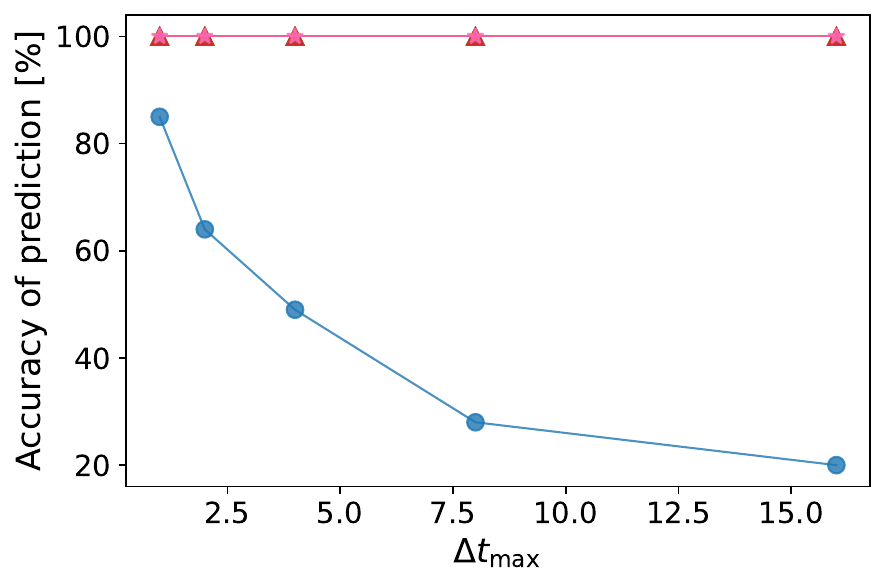}
    \end{subfigure}
    \begin{subfigure}{0.48\textwidth}
        \centering
        \includegraphics[width=\textwidth]{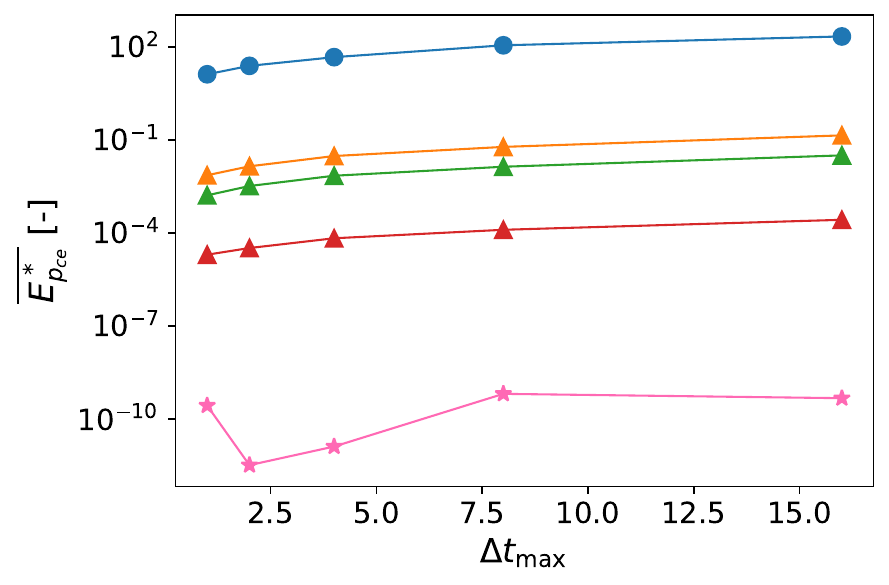}
    \end{subfigure}
    \caption{\textbf{Comparison of prediction accuracy between \FIN, \FIR, and \FIT.} Left: probability (accuracy) of correctly predicting invasion pattern for the different methods over maximum time-step size. Right: invasion event prediction errors, \cref{eq:errorinvasions}, over maximum time-step size.}
	\label{fig:bifurcated_benchmark_case}
\end{figure}

We generate 100 realizations of the network shown in \cref{fig:set_up_bifuracted}, with randomly generated throat radii for each realization. 
We then vary the maximum allowed time-step size to analyze the convergence behavior. Here, we choose the maximum time-step size, since $\overline{\Delta t}$ may vary across the 100 realizations (see \cref{re:newtontimesteps}). \cref{fig:bifurcated_benchmark_case} (left) shows the probability of correctly predicting the invasion sequence for different time-step sizes. \cref{fig:bifurcated_benchmark_case} (right) shows the prediction error, \cref{eq:errorinvasions}, averaged over all realizations for which the invasion sequence was correctly predicted.  

Similar to the previous analysis with the one-dimensional network, the \FIT scheme exactly matches the entry pressure for invasion events such that there are no prediction errors, as shown in \cref{fig:bifurcated_benchmark_case} (right). All other schemes do have prediction errors within time steps, where invasion events occur. By reducing the regularization interval $\delta$, these errors can be reduced for the \FIR schemes. 
The prediction errors of the \FIN scheme are approximately three orders of magnitude larger than those of \FIR  with $\delta = 0.4$.

The prediction probability of the \FIN scheme decreases significantly for larger time-step sizes, as shown in \cref{fig:bifurcated_benchmark_case} (left), and never exceeds $90\%$.
All other schemes consistently predict the invasion sequences correctly, regardless of the maximum time-step size. The large errors of the \FIN scheme are caused by incorrect decisions made during Newton iterations based on unconverged solutions, as previously discussed.

\subsection{Two-dimensional fully-connected network case} 
\begin{figure}[h!]
    \centering
    \includegraphics[width=\linewidth]{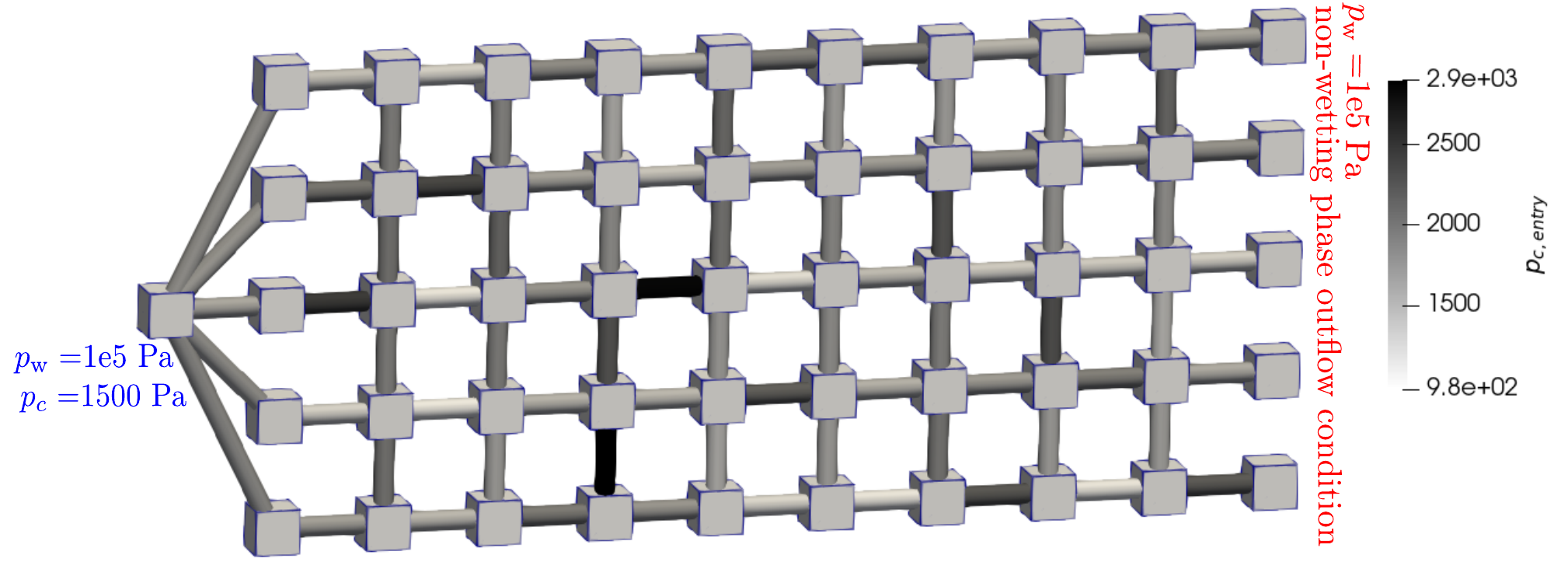}
    \caption{\textbf{Setup for the 2D drainage case.} Non-wetting phase enters the domain through the inlet pore, where a capillary pressure of $p_{ce, \mathrm{inlet}} = 1.5 \times 10^2$ \si{\pascal} is set.}
    \label{fig:3d_drainage_setup}
\end{figure}

In this test case, a two-dimensional drainage problem characterized by a global capillary pressure difference imposed between the inlet and outlet is considered for analysis. At the outlet, the non-wetting and wetting phase can flow out and the wetting-phase pressure is fixed $p_\mathrm{w, outlet} = 10^5$ \si{\pascal}. At the inlet, the wetting phase pressure is prescribed as $p_\mathrm{w, inlet} = 10^5$ \si{\pascal}, and the non-wetting phase pressure is computed accordingly by $p_\mathrm{n, inlet} = p_\mathrm{w, inlet} + p_{ce, \mathrm{inlet}}$.  To construct the network, we generate a lattice network consisting of 50 cubic pore bodies arranged in 10 columns and 5 rows, which are connected by throats with randomly distributed radii sampled from a Gauss normal distribution. Then, this lattice network is connected to an inlet pore body at one side via throats with similar entry capillary pressure of $p_{ce, \mathrm{inlet}} = 1.5 \times 10^2$ \si{\pascal}. The problem setup is illustrated in \cref{fig:3d_drainage_setup}. To investigate the behavior of the presented schemes, 10 realizations, i.e., $S1$ to $S10$, with randomly generated throat radii are considered in the following analysis.

Instead of the prediction error \cref{eq:errorinvasions}, we use
\begin{equation}
\label{eq:errorinvasionsmax}
    E^{*,\max}_{p_{ce}} := \max_{\lbrace t^{k},ij \rbrace \in \mathcal{T}^*_N}  \frac{| \dpthroat^{k} |}{p_{ce,ij}}.
\end{equation}
The reason for this is that, using the mean error \cref{eq:errorinvasions} for this more complex test case can not clearly show the influence of certain local effects that occur in specific throats of the network, as can be seen later in the discussion at the end of this section.

\begin{figure}[h!]
    \centering
    \includegraphics[width=0.6\columnwidth]{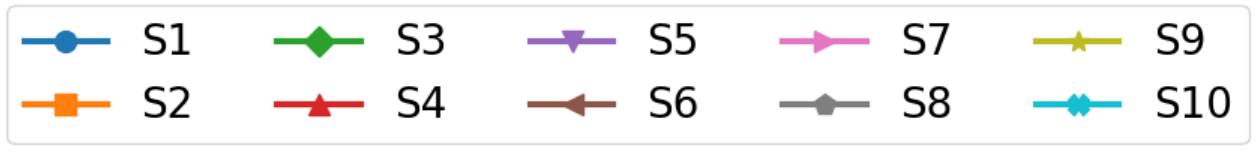}

    \begin{subfigure}{\columnwidth}
        \includegraphics[width=\columnwidth]{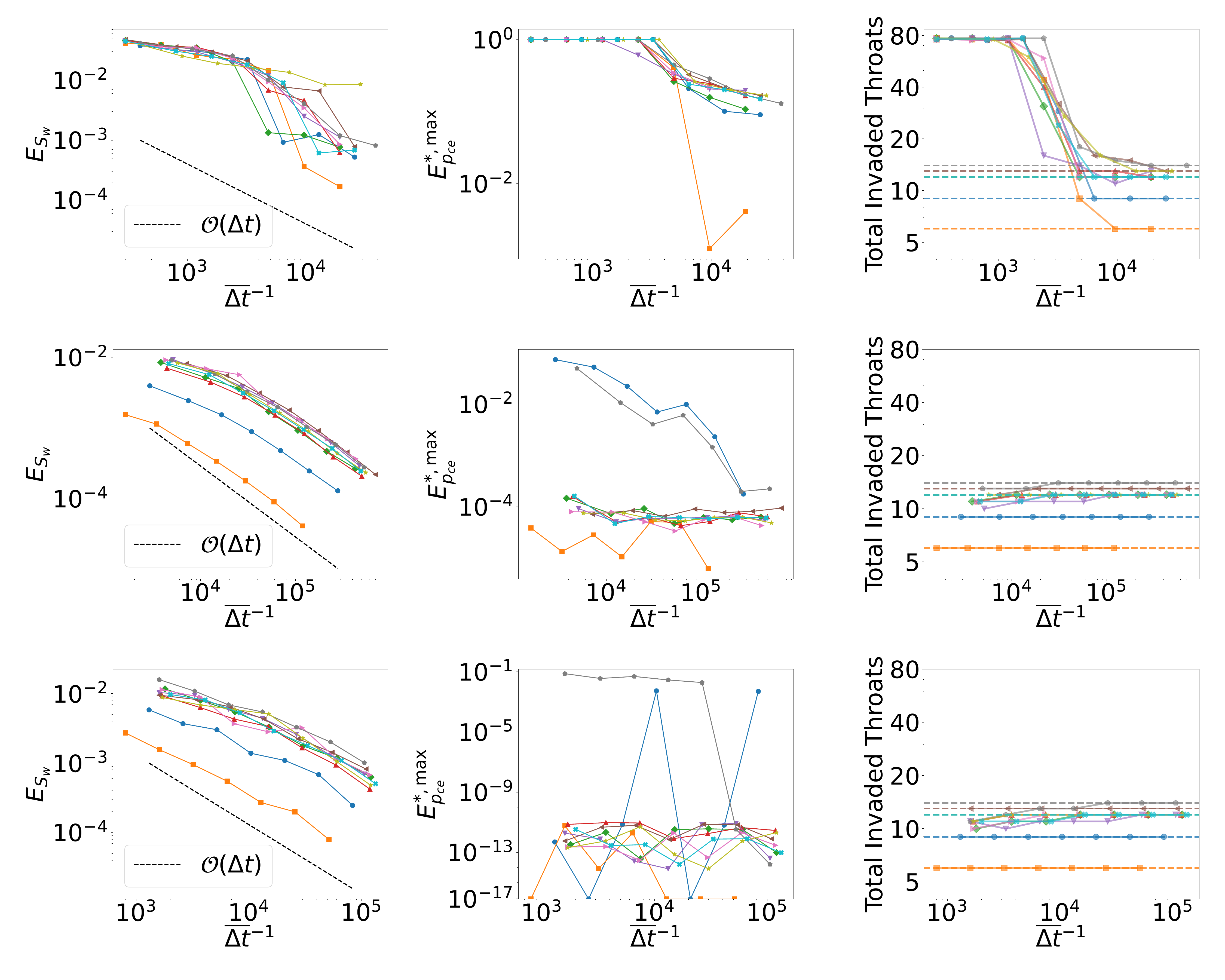}
    \end{subfigure}

    \caption{\textbf{Temporal convergence.} The first column presents the $L^2$-errors with a reference line of first-order convergence (black dashed line). The second column shows $E^{*,\max}_{p_{ce}}$, \cref{eq:errorinvasionsmax}, over $\overline{\Delta t}^{-1}$. The third column depicts the total number of invaded throats, where the dashed lines correspond to the number of the reference solutions. The first row are the results of the \FIN, the second row of \FIR with $\delta = 10^{-3}$, and the third row of \FIT schemes. }
    \label{fig:temporal-convergence}
\end{figure}

\cref{fig:temporal-convergence} shows the $L^2$-errors (left column), the errors related to invasion events (middle column), i.e. prediction errors, and the total number of invaded throats at the end of the simulations (right column). The first row presents the results of the \FIN, the second row of \FIR (with $\delta = 10^{-3}$), and the third row of \FIT schemes.
Examining the $L^2$-errors, it can be observed that the \FIR and \FIT schemes demonstrate convergence rates that are close to first order. In contrast to the \FIN scheme, where the errors fluctuate and there is no clear convergence trend. 
This is caused by inaccuracies in predicting invasion events, as shown in the second column, where the \FIN scheme, particularly for small $\overline{\Delta t}$, exhibits larger prediction errors compared to the other schemes. As previously discussed, these errors significantly affect the overall $L^2$-errors $E_{S_w}$, and convergence can generally only be expected if the method is able to control the prediction errors. This is an inherent ability of \FIT, where we imposed the additional throat constraints and also for \FIR, when small regularization intervals $\delta$ are chosen. 
For two cases, the prediction errors do not vanish (solver tolerance) for the \FIT scheme, which will be discussed at the end of this section. 

The right column of \cref{fig:temporal-convergence} shows that the \FIR and \FIT schemes predict the total number of invaded throats with high accuracy, even for larger time-step sizes. Please note that for larger time steps, multiple throats can be invaded within a single time step, especially if invasion time points are close together. When not resolving each invasion separately, deviations from the exact invasion sequence of the reference solution may occur. However, the overall behavior is still well captured, which is in contrast to the behavior of the \FIN scheme, where it predicts a completely wrong invasion pattern for larger time-step sizes. 

\todo[inline]{We could plot some averaged results to show that such quantities still match, even thought local patterns are completely wrong. }

\begin{figure}[h!]
    \centering
    \includegraphics[width=\linewidth]{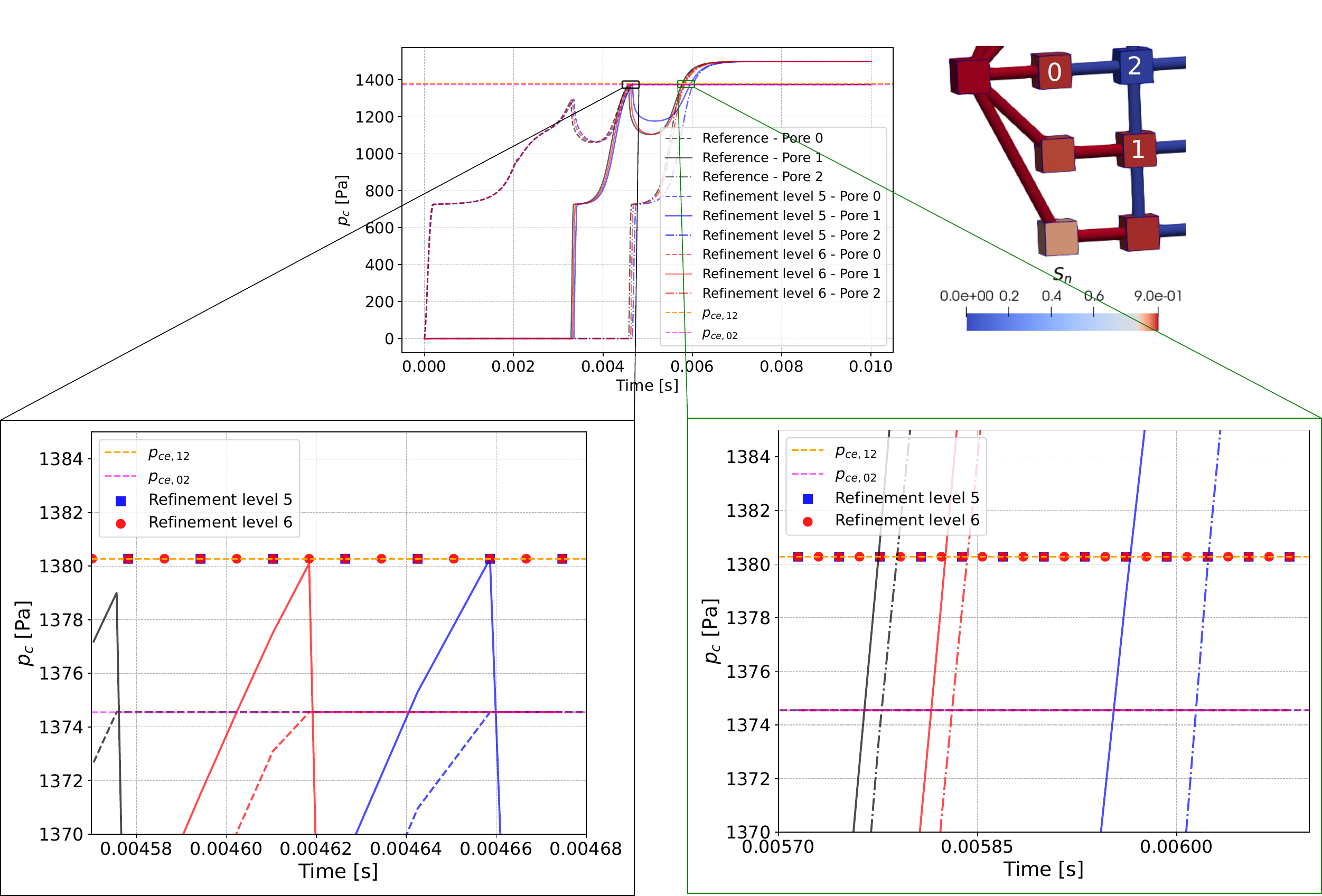}
    \caption{\textbf{Capillary pressure over time} plotted for three pore bodies of the network for different time refinement levels. The results are shown for the \FIT scheme. The lower plots show a zoom of the overall time evolution. The results correspond to scenario $S1$ of \cref{fig:temporal-convergence}.}
    \label{fig:enter-label}
\end{figure}

In the last part of this section, we discuss why the prediction errors of the \FIT scheme (\cref{fig:temporal-convergence}, bottom middle) are not always zero (solver tolerance). Examining the additional throat constraint \cref{eq:timedisc_throat_md}, the only reason that $\dpthroat \not = 0$ for an invasion event is that $\throatvar_{ij} = 1$ (see \cref{eq:state_invadedbeginning}). In the previous test cases, this did not occur because the smaller the coordination number (i.e. $|\calN_i|$), the lower is the probability of encountering such a behavior, as for the one-dimensional network and also for most of the pore bodies in the bifurcating network we had a coordination number of two. Otherwise, for each invasion event, $\dpthroat = 0$.
In \cref{fig:enter-label}, the time evolution of the capillary pressure solution of three pore bodies is plotted for the \FIT scheme on refinement levels 5 and 6 and compared to the reference solution (calculated for even smaller time-step sizes). It can be observed that pore bodies 0 and 1 are earlier invaded by the non-wetting fluid. At approximately $0.0046 \pm \delta t$ seconds, the capillary pressures of both of these pore bodies are close to the throat entry capillary pressures, required to invade the throat connecting them to pore body 2, see \cref{fig:enter-label} lower left plot. For the calculations on refinement level 5, both throats become invaded within one time step. However, after time-step refinement, invasion occurs only in throat $\lbrace 0,2 \rbrace$.
As mentioned earlier, such events are of non-local effects, resulting in water fluxes into pore body 1 due to water displacement in pore body 2.

So far, $\dpthroatidx{12} \leq 0$, which does not yet explain why $\throatvar_{12} = 1$. However, this can be explained by looking at the lower right plot, which shows the capillary pressure evolution at times $0.0059 \pm \delta t$ seconds. The non-wetting phase influx from pore 0 displaced the wetting phase from pore 2, causing its capillary pressure to increase to nearly the capillary pressure of pore 1. As a result, after the throat between pores 1 and 2 becomes invaded, there is no relaxation time, since the resistance, due to the already high capillary pressure of pore 2, is too high. Therefore, even though  $\throatvar_{12} = 1$, the non-wetting flux flowing from pore 1 to 2 is not large enough (small capillary pressure differences) to maintain $\dpthroatidx{12} \leq 0$.
This explains the fluctuating behavior of the results for $S1$ shown in \cref{fig:temporal-convergence} middle-bottom plot. Similar behavior is also observed for scenario $S8$. 
This means that for networks with high connectivity, i.e, high average coordination number, competing invasion events and non-local effects may lead to such behavior. 
In such cases, one should not refer to prediction error if $\dpthroatidx{12} > 0$ and $\dpthroatidx{12} = 1$. Therefore, prediction errors should only be related to cases where $\dpthroatidx{12} > 0$ and $\dpthroatidx{12} \in (0,1)$, which is prevented by the \FIT scheme.
In general, the plots shown in \cref{fig:enter-label} again highlight the good convergence behavior of the \FIT scheme, even though the local invasion pattern differs between refinement levels 5 and 6.

\paragraph{Summary:}
\begin{itemize}
    \item The overall $L^2$-error is strongly influenced by the prediction errors.
    \item Only if errors for each local event (here invasion) are controlled, overall convergence can be expected. 
    \item The new schemes are accurate also for larger time-step sizes and show the expected convergence behavior. 
    \item The \FIN scheme becomes unreliable (completely wrong invasion pattern) for large time steps due to incorrect decisions made during the Newton iterations, based on unconverged solutions.
\end{itemize}
The behavior where single local time events (such as invasions or snap-off) may dominate the overall error is specific to the type of models considered here, due to the discontinuities in conductivities, which distinguishes them from classical models having no discontinuities.

\todo[inline]{We always talk about phases but don't really specify details about the non-wetting phase? Is it air? I guess it would be interesting to at least get a feeling about the properties of the non-wetting fluid}

\section{Summary and Outlook}
In this manuscript, we have discussed in detail over- and underestimation of fluxes during invasion (snap-off) events when using a fully-implicit Euler time discretization scheme. This has been highlighted by considering the exact time-averaged fluxes, where a scaling factor $\throatvar$ has been introduced. It has been shown that $\throatvar \rightarrow 1$ if the event (invasion or snap-off) happens at the beginning of the time step and that  $\throatvar \rightarrow 0$ if it happens at the end. This, for example, also means that the fluxes $Q \rightarrow 0$ if $\throatvar \rightarrow 0$. Such behavior is not reflected by classical fully-implicit flux approximation schemes, thus, leading to large flux approximation errors and related event prediction errors. 
To overcome these issues, we have introduced, see \cref{sec:GenFluxFunction}, the concept of generalized flux function which depend on the scaling factor $\throatvar$. In \cref{subsec:fluxreg}, this factor has then been expressed by introducing regularized Heaviside functions. In \cref{subsec:mdtheta}, it has been kept as an additional throat variable and an additional equation has been formulated to express the corresponding throat state. This allows exactly capturing the invasion criterion within the time step where such events happen. Furthermore, the different states are discussed in detail. This new approach leads to a coupled system between mass balance equations for pore bodies and state equations for throats. It should be mentioned that solving such coupled systems of equations is naturally more numerically demanding, requiring the development of more suitable linear solvers. As previously mentioned, trivial states (no invasion) could also be eliminated, further simplifying the coupled system. However, addressing these aspects is beyond the scope of this paper and will be considered in future work. 

To analyze the temporal convergence behavior of the new schemes and to compare it with an established one, three test cases with increasing complexity have been considered. 
The first one, a one-dimensional homogeneous network, allows constructing an analytical solution. 
In the second test case, a bifurcating network has been considered, where throat radii are randomly generated leading to random capillary entry pressures. Further extensions have been made in the third test case, where a larger fully-connected network has been considered. For all test cases, the overall $L^2$-errors and the related temporal convergence have been analyzed. In addition, prediction errors related to invasion events have been investigated and compared for the different schemes. It has been demonstrated that the \FIT scheme accurately matches the capillary entry pressures during invasion events and that the new schemes exhibit smaller errors and a clear convergence behavior compared to the \FIN scheme. 
Regardless of the time-step sizes used, no instabilities have been observed, indicating the expected unconditional stability of the presented fully-implicit schemes.
Furthermore, the importance of controlling prediction errors to obtain overall temporal convergence has been discussed in detail. Such behavior is specific to the pore-network models considered in this work, which exhibit discontinuities in conductivities during invasion events. 
To control prediction errors and improve the efficiency of the schemes, applying local time-stepping adaptation is beneficial for resolving invasion (snap-off) events with smaller time-step sizes. Such a strategy has been applied for calculating the reference solutions. However, in this work, we aimed to focus exclusively on the convergence behavior of the presented schemes without complicating the discussion by incorporating such strategies. This is also the focus of ongoing work, where accurate time predictors are being investigated. 

In future work, we also aim to explore the mathematical theory behind the presented models, which could serve as the basis for proving temporal convergence.



\bmhead{Supplementary information}
\bmhead{Acknowledgments}
We acknowledge funding by the German Research Foundation (DFG), within the Collaborative Research Center on Interface-Driven Multi-Field Processes in Porous Media (SFB 1313, Project No. 327154368). 
We also acknowledge Special Research Fund (BOF) of Hasselt University  (Project BOF22KV03) and
Research Foundation-Flanders (FWO), Belgium (Project G0A9A25N).
\bmhead{Funding}
German Research Foundation (DFG), SFB 1313, Project No. 327154368. \\
Special Research Fund (BOF) of Hasselt University, Project BOF22KV03 \\
Research Foundation-Flanders (FWO), Belgium, Project G0A9A25N










\begin{appendices}






\end{appendices}


\bibliography{refs}


\begin{thebibliography}{22}
\ifx \bisbn   \undefined \def \bisbn  #1{ISBN #1}\fi
\ifx \binits  \undefined \def \binits#1{#1}\fi
\ifx \bauthor  \undefined \def \bauthor#1{#1}\fi
\ifx \batitle  \undefined \def \batitle#1{#1}\fi
\ifx \bjtitle  \undefined \def \bjtitle#1{#1}\fi
\ifx \bvolume  \undefined \def \bvolume#1{\textbf{#1}}\fi
\ifx \byear  \undefined \def \byear#1{#1}\fi
\ifx \bissue  \undefined \def \bissue#1{#1}\fi
\ifx \bfpage  \undefined \def \bfpage#1{#1}\fi
\ifx \blpage  \undefined \def \blpage #1{#1}\fi
\ifx \burl  \undefined \def \burl#1{\textsf{#1}}\fi
\ifx \doiurl  \undefined \def \doiurl#1{\url{https://doi.org/#1}}\fi
\ifx \betal  \undefined \def \betal{\textit{et al.}}\fi
\ifx \binstitute  \undefined \def \binstitute#1{#1}\fi
\ifx \binstitutionaled  \undefined \def \binstitutionaled#1{#1}\fi
\ifx \bctitle  \undefined \def \bctitle#1{#1}\fi
\ifx \beditor  \undefined \def \beditor#1{#1}\fi
\ifx \bpublisher  \undefined \def \bpublisher#1{#1}\fi
\ifx \bbtitle  \undefined \def \bbtitle#1{#1}\fi
\ifx \bedition  \undefined \def \bedition#1{#1}\fi
\ifx \bseriesno  \undefined \def \bseriesno#1{#1}\fi
\ifx \blocation  \undefined \def \blocation#1{#1}\fi
\ifx \bsertitle  \undefined \def \bsertitle#1{#1}\fi
\ifx \bsnm \undefined \def \bsnm#1{#1}\fi
\ifx \bsuffix \undefined \def \bsuffix#1{#1}\fi
\ifx \bparticle \undefined \def \bparticle#1{#1}\fi
\ifx \barticle \undefined \def \barticle#1{#1}\fi
\bibcommenthead
\ifx \bconfdate \undefined \def \bconfdate #1{#1}\fi
\ifx \botherref \undefined \def \botherref #1{#1}\fi
\ifx \url \undefined \def \url#1{\textsf{#1}}\fi
\ifx \bchapter \undefined \def \bchapter#1{#1}\fi
\ifx \bbook \undefined \def \bbook#1{#1}\fi
\ifx \bcomment \undefined \def \bcomment#1{#1}\fi
\ifx \oauthor \undefined \def \oauthor#1{#1}\fi
\ifx \citeauthoryear \undefined \def \citeauthoryear#1{#1}\fi
\ifx \endbibitem  \undefined \def \endbibitem {}\fi
\ifx \bconflocation  \undefined \def \bconflocation#1{#1}\fi
\ifx \arxivurl  \undefined \def \arxivurl#1{\textsf{#1}}\fi
\csname PreBibitemsHook\endcsname

\bibitem[\protect\citeauthoryear{Blunt}{2001}]{blunt2001flow}
\begin{barticle}
\bauthor{\bsnm{Blunt}, \binits{M.J.}}:
\batitle{Flow in porous media—pore-network models and multiphase flow}.
\bjtitle{Current opinion in colloid \& interface science}
\bvolume{6}(\bissue{3}),
\bfpage{197}--\blpage{207}
(\byear{2001})
\end{barticle}
\endbibitem

\bibitem[\protect\citeauthoryear{Raoof et~al.}{2013}]{raoof2013poreflow}
\begin{barticle}
\bauthor{\bsnm{Raoof}, \binits{A.}},
\bauthor{\bsnm{Nick}, \binits{H.M.}},
\bauthor{\bsnm{Hassanizadeh}, \binits{S.M.}},
\bauthor{\bsnm{Spiers}, \binits{C.}}:
\batitle{Poreflow: A complex pore-network model for simulation of reactive
  transport in variably saturated porous media}.
\bjtitle{Computers \& Geosciences}
\bvolume{61},
\bfpage{160}--\blpage{174}
(\byear{2013})
\end{barticle}
\endbibitem

\bibitem[\protect\citeauthoryear{Chen et~al.}{2020}]{chen2020fully}
\begin{barticle}
\bauthor{\bsnm{Chen}, \binits{S.}},
\bauthor{\bsnm{Qin}, \binits{C.}},
\bauthor{\bsnm{Guo}, \binits{B.}}:
\batitle{Fully implicit dynamic pore-network modeling of two-phase flow and
  phase change in porous media}.
\bjtitle{Water Resources Research}
\bvolume{56}(\bissue{11}),
\bfpage{2020}--\blpage{028510}
(\byear{2020})
\end{barticle}
\endbibitem

\bibitem[\protect\citeauthoryear{Veyskarami
  et~al.}{2016}]{veyskarami2016modeling}
\begin{barticle}
\bauthor{\bsnm{Veyskarami}, \binits{M.}},
\bauthor{\bsnm{Hassani}, \binits{A.H.}},
\bauthor{\bsnm{Ghazanfari}, \binits{M.H.}}:
\batitle{Modeling of non-darcy flow through anisotropic porous media: Role of
  pore space profiles}.
\bjtitle{Chemical Engineering Science}
\bvolume{151},
\bfpage{93}--\blpage{104}
(\byear{2016})
\end{barticle}
\endbibitem

\bibitem[\protect\citeauthoryear{Celia et~al.}{1995}]{celia1995recent}
\begin{barticle}
\bauthor{\bsnm{Celia}, \binits{M.A.}},
\bauthor{\bsnm{Reeves}, \binits{P.C.}},
\bauthor{\bsnm{Ferrand}, \binits{L.A.}}:
\batitle{Recent advances in pore scale models for multiphase flow in porous
  media}.
\bjtitle{Reviews of Geophysics}
\bvolume{33}(\bissue{S2}),
\bfpage{1049}--\blpage{1057}
(\byear{1995})
\end{barticle}
\endbibitem

\bibitem[\protect\citeauthoryear{Koplik and Lasseter}{1985}]{koplik1985two}
\begin{barticle}
\bauthor{\bsnm{Koplik}, \binits{J.}},
\bauthor{\bsnm{Lasseter}, \binits{T.}}:
\batitle{Two-phase flow in random network models of porous media}.
\bjtitle{Society of Petroleum Engineers Journal}
\bvolume{25}(\bissue{01}),
\bfpage{89}--\blpage{100}
(\byear{1985})
\end{barticle}
\endbibitem

\bibitem[\protect\citeauthoryear{Thompson}{2002}]{thompson2002pore}
\begin{barticle}
\bauthor{\bsnm{Thompson}, \binits{K.E.}}:
\batitle{Pore-scale modeling of fluid transport in disordered fibrous
  materials}.
\bjtitle{AIChE journal}
\bvolume{48}(\bissue{7}),
\bfpage{1369}--\blpage{1389}
(\byear{2002})
\end{barticle}
\endbibitem

\bibitem[\protect\citeauthoryear{Joekar-Niasar et~al.}{2010}]{joekar2010non}
\begin{barticle}
\bauthor{\bsnm{Joekar-Niasar}, \binits{V.}},
\bauthor{\bsnm{Hassanizadeh}, \binits{S.M.}},
\bauthor{\bsnm{Dahle}, \binits{H.}}:
\batitle{Non-equilibrium effects in capillarity and interfacial area in
  two-phase flow: dynamic pore-network modelling}.
\bjtitle{Journal of fluid mechanics}
\bvolume{655},
\bfpage{38}--\blpage{71}
(\byear{2010})
\end{barticle}
\endbibitem

\bibitem[\protect\citeauthoryear{Weishaupt and Helmig}{2021}]{Weishaupt2021}
\begin{barticle}
\bauthor{\bsnm{Weishaupt}, \binits{K.}},
\bauthor{\bsnm{Helmig}, \binits{R.}}:
\batitle{A dynamic and fully implicit non-isothermal, two-phase, two-component
  pore-network model coupled to single-phase free flow for the pore-scale
  description of evaporation processes}.
\bjtitle{Water Resources Research}
\bvolume{57}(\bissue{4}),
\bfpage{2020}--\blpage{028772}
(\byear{2021})
\doiurl{10.1029/2020WR028772}
\end{barticle}
\endbibitem

\bibitem[\protect\citeauthoryear{An et~al.}{2020}]{an2020transition}
\begin{barticle}
\bauthor{\bsnm{An}, \binits{S.}},
\bauthor{\bsnm{Erfani}, \binits{H.}},
\bauthor{\bsnm{Godinez-Brizuela}, \binits{O.E.}},
\bauthor{\bsnm{Niasar}, \binits{V.}}:
\batitle{Transition from viscous fingering to capillary fingering: Application
  of gpu-based fully implicit dynamic pore network modeling}.
\bjtitle{Water Resources Research}
\bvolume{56}(\bissue{12}),
\bfpage{2020}--\blpage{028149}
(\byear{2020})
\end{barticle}
\endbibitem

\bibitem[\protect\citeauthoryear{Wu et~al.}{2024}]{Wu2024}
\begin{barticle}
\bauthor{\bsnm{Wu}, \binits{H.}},
\bauthor{\bsnm{Veyskarami}, \binits{M.}},
\bauthor{\bsnm{Schneider}, \binits{M.}},
\bauthor{\bsnm{Helmig}, \binits{R.}}:
\batitle{A new fully implicit two-phase pore-network model by utilizing
  regularization strategies}.
\bjtitle{Transport in Porous Media}
\bvolume{151}(\bissue{1}),
\bfpage{1}--\blpage{26}
(\byear{2024})
\doiurl{10.1007/s11242-023-02031-2}
\end{barticle}
\endbibitem

\bibitem[\protect\citeauthoryear{Weishaupt}{2020}]{weishaupt2020model}
\begin{bbook}
\bauthor{\bsnm{Weishaupt}, \binits{K.}}:
\bbtitle{Model Concepts for Coupling Free Flow with Porous Medium Flow at the
  Pore-network Scale: from Single-phase Flow to Compositional Non-isothermal
  Two-phase Flow}.
\bpublisher{Stuttgart: Eigenverlag des Instituts f{\"u}r Wasser-und
  Umweltsystemmodellierung der Universität Stuttgart}, \blocation{???}
(\byear{2020})
\end{bbook}
\endbibitem

\bibitem[\protect\citeauthoryear{Dias and Payatakes}{1986a}]{dias1986network}
\begin{barticle}
\bauthor{\bsnm{Dias}, \binits{M.M.}},
\bauthor{\bsnm{Payatakes}, \binits{A.C.}}:
\batitle{Network models for two-phase flow in porous media part 1. immiscible
  microdisplacement of non-wetting fluids}.
\bjtitle{Journal of Fluid Mechanics}
\bvolume{164},
\bfpage{305}--\blpage{336}
(\byear{1986})
\end{barticle}
\endbibitem

\bibitem[\protect\citeauthoryear{Dias and Payatakes}{1986b}]{dias1986network2}
\begin{barticle}
\bauthor{\bsnm{Dias}, \binits{M.M.}},
\bauthor{\bsnm{Payatakes}, \binits{A.C.}}:
\batitle{Network models for two-phase flow in porous media part 2. motion of
  oil ganglia}.
\bjtitle{Journal of Fluid Mechanics}
\bvolume{164},
\bfpage{337}--\blpage{358}
(\byear{1986})
\end{barticle}
\endbibitem

\bibitem[\protect\citeauthoryear{Al-Gharbi and Blunt}{2005}]{al2005dynamic}
\begin{barticle}
\bauthor{\bsnm{Al-Gharbi}, \binits{M.S.}},
\bauthor{\bsnm{Blunt}, \binits{M.J.}}:
\batitle{Dynamic network modeling of two-phase drainage in porous media}.
\bjtitle{Physical Review E—Statistical, Nonlinear, and Soft Matter Physics}
\bvolume{71}(\bissue{1}),
\bfpage{016308}
(\byear{2005})
\end{barticle}
\endbibitem

\bibitem[\protect\citeauthoryear{Joekar~Niasar
  et~al.}{2009}]{joekar2009simulating}
\begin{botherref}
\oauthor{\bsnm{Joekar~Niasar}, \binits{V.}},
\oauthor{\bsnm{Hassanizadeh}, \binits{S.}},
\oauthor{\bsnm{Pyrak-Nolte}, \binits{L.}},
\oauthor{\bsnm{Berentsen}, \binits{C.}}:
Simulating drainage and imbibition experiments in a high-porosity micromodel
  using an unstructured pore network model.
Water resources research
\textbf{45}(2)
(2009)
\end{botherref}
\endbibitem

\bibitem[\protect\citeauthoryear{Berg et~al.}{2013}]{berg2013real}
\begin{barticle}
\bauthor{\bsnm{Berg}, \binits{S.}},
\bauthor{\bsnm{Ott}, \binits{H.}},
\bauthor{\bsnm{Klapp}, \binits{S.A.}},
\bauthor{\bsnm{Schwing}, \binits{A.}},
\bauthor{\bsnm{Neiteler}, \binits{R.}},
\bauthor{\bsnm{Brussee}, \binits{N.}},
\bauthor{\bsnm{Makurat}, \binits{A.}},
\bauthor{\bsnm{Leu}, \binits{L.}},
\bauthor{\bsnm{Enzmann}, \binits{F.}},
\bauthor{\bsnm{Schwarz}, \binits{J.-O.}}, \betal:
\batitle{Real-time 3d imaging of haines jumps in porous media flow}.
\bjtitle{Proceedings of the National Academy of Sciences}
\bvolume{110}(\bissue{10}),
\bfpage{3755}--\blpage{3759}
(\byear{2013})
\end{barticle}
\endbibitem

\bibitem[\protect\citeauthoryear{Sun and Santamarina}{2019}]{sun2019haines}
\begin{barticle}
\bauthor{\bsnm{Sun}, \binits{Z.}},
\bauthor{\bsnm{Santamarina}, \binits{J.C.}}:
\batitle{Haines jumps: Pore scale mechanisms}.
\bjtitle{Physical review E}
\bvolume{100}(\bissue{2}),
\bfpage{023115}
(\byear{2019})
\end{barticle}
\endbibitem

\bibitem[\protect\citeauthoryear{Gjennestad
  et~al.}{2018}]{gjennestad2018stable}
\begin{barticle}
\bauthor{\bsnm{Gjennestad}, \binits{M.A.}},
\bauthor{\bsnm{Vassvik}, \binits{M.}},
\bauthor{\bsnm{Kjelstrup}, \binits{S.}},
\bauthor{\bsnm{Hansen}, \binits{A.}}:
\batitle{Stable and efficient time integration of a dynamic pore network model
  for two-phase flow in porous media}.
\bjtitle{Frontiers in Physics}
\bvolume{6},
\bfpage{56}
(\byear{2018})
\end{barticle}
\endbibitem

\bibitem[\protect\citeauthoryear{Aker et~al.}{1998}]{aker1998}
\begin{barticle}
\bauthor{\bsnm{Aker}, \binits{E.}},
\bauthor{\bsnm{J{\O}rgen~M{\AA}l{\O}y}, \binits{K.}},
\bauthor{\bsnm{Hansen}, \binits{A.}},
\bauthor{\bsnm{Batrouni}, \binits{G.G.}}:
\batitle{A two-dimensional network simulator for two-phase flow in porous
  media}.
\bjtitle{Transport in Porous Media}
\bvolume{32}(\bissue{2}),
\bfpage{163}--\blpage{186}
(\byear{1998})
\doiurl{10.1023/A:1006510106194}
\end{barticle}
\endbibitem

\bibitem[\protect\citeauthoryear{Akkerman}{2017}]{Akkerman2017}
\begin{barticle}
\bauthor{\bsnm{Akkerman}, \binits{I.}}:
\batitle{Monotone level-sets on arbitrary meshes without redistancing}.
\bjtitle{Computers \& Fluids}
\bvolume{146},
\bfpage{74}--\blpage{85}
(\byear{2017})
\doiurl{10.1016/j.compfluid.2017.01.007}
\end{barticle}
\endbibitem

\bibitem[\protect\citeauthoryear{Armstrong and Berg}{2013}]{armstrong2013}
\begin{barticle}
\bauthor{\bsnm{Armstrong}, \binits{R.T.}},
\bauthor{\bsnm{Berg}, \binits{S.}}:
\batitle{Interfacial velocities and capillary pressure gradients during haines
  jumps}.
\bjtitle{Phys. Rev. E}
\bvolume{88},
\bfpage{043010}
(\byear{2013})
\doiurl{10.1103/PhysRevE.88.043010}
\end{barticle}
\endbibitem

\end{thebibliography}

\end{document}